\documentclass[10pt]{amsart}
\usepackage{indentfirst}
\usepackage{amsmath,amssymb,latexsym,esint,cite,mathrsfs}
\usepackage{verbatim,wasysym}
\usepackage[left=2.6cm,right=2.6cm,top=2.9cm,bottom=2.7cm]{geometry}
\usepackage{tikz,enumitem,graphicx, subfig, microtype, color}
\usepackage{epic,eepic}
\usepackage[colorlinks=true,urlcolor=blue, citecolor=red,linkcolor=blue,
linktocpage,pdfpagelabels, bookmarksnumbered,bookmarksopen]{hyperref}
\usepackage[hyperpageref]{backref}
\usepackage[english]{babel}
\numberwithin{equation}{section}
\newtheorem{thm}{Theorem}[section]
\newtheorem{lem}[thm]{Lemma}

\newtheorem{Rem}[thm]{Remark}

\renewcommand{\i}{{\mathrm i}}

\begin{document}
\title[Helmholtz equation]{Complex and Real valued solutions for Fractional Helmholtz equation}
\author{ Zifei Shen$^1$ and Shuijin Zhang$^{1,\ast}$}

.

\begin{abstract}
In this paper, we are concerned with the limiting absorption principle for the fractional Helmholtz equation
\begin{equation}\label{main system}
(-\Delta)^{s} u-\lambda u=f(x,u), ~~\mathrm{in}~~\mathbb{R}^{n},
\end{equation}
where $n\geq3$, $0<\lambda<+\infty$ and $\frac{n}{n+1}<s<\frac{n}{2}$ are two real parameters. By establishing the boundedness estimate for the resolvent of fractional Helmholtz operator, we obtain the nontrivial $L^{q}(\mathbb{R}^{n})$ complex valued solutions for (\ref{main system}). By setting up a dual variational framework, we also obtain the real valued solutions for (\ref{main system}) via a non-vanishing principle.
\end{abstract}

\renewcommand{\thefootnote}{\fnsymbol{footnote}}
\footnotetext{\hspace*{-5mm}
\renewcommand{\arraystretch}{1}
\begin{tabular}{@{}r@{}p{12cm}@{}}
$^\ast$& ~Corresponding author, e-mail:shuijinzhang@zjnu.edu.cn;
\end{tabular}}

\renewcommand{\thefootnote}{\fnsymbol{footnote}}
\footnotetext{\hspace*{-5mm}
\renewcommand{\arraystretch}{1}
\begin{tabular}{@{}r@{}p{12cm}@{}}
~~~&~Zifei Shen: szf@zjnu.edu.cn
\end{tabular}}

\renewcommand{\thefootnote}{\fnsymbol{footnote}}
\footnotetext{\hspace*{-5mm}
\renewcommand{\arraystretch}{1}
\begin{tabular}{@{}r@{}p{12cm}@{}}
$^1$& ~Department of Mathematics, Zhejiang Normal University, 321000, Jinhua, China.
\end{tabular}}
~\\
   
\maketitle
\setlength{\parindent}{2em}
\section{Introduction and Main Results}
In this paper, we are concerned with the nontrivial $L^{q}(\mathbb{R}^{n})$ complex valued solutions for the fractional Helmholtz equation
\begin{equation}\label{main system 1}
(-\Delta)^{s} u-\lambda u=f(x,u), ~~\mathrm{in}~~\mathbb{R}^{n},
\end{equation}
where $n\geq3$, $0<\lambda<+\infty$ and $\frac{n}{n+1}<s<\frac{n}{2}$ are two real parameters.

As $s=1$, the Helmholtz equation
\begin{equation}\label{single}
-\Delta u-\lambda u=f(x,u), ~~~\mathrm{in}~\mathbb{R}^{n},
\end{equation}
arises frequently in scattering theory and quantum mechanics.
Since the parameter $\lambda>0$ is contained in the essential spectrum
of negative Laplacian, the Helmholtz equation performs a strongly indefinite structure.
If $f(x,u)$ is assumed to be a local nonlinearity,
that is $f(x,u)=0$ in $\mathbb{R}^{n}\setminus B_{R}$,
one can investigate the Helmholtz equation in a bounded domain based on a Dirichlet to Neumann map,
and by a Linking argument, one can obtain the real valued solutions $u\in H^{1}_{loc}(\mathbb{R}^{n})$
for (\ref{single}), see Ev\'{e}quoz and Weth \cite{Evequoz2014}.
However, if $f(x,u)$ is a physically source term without compact support in $\mathbb{R}^{n}$, $H^{1}(\mathbb{R}^{n})$ may be not a proper space to find solutions. More precisely, since the Helmholtz operator is defined on the whole space and hence the oscillating solutions for homogeneous Helmholtz equation with slow decay which, in general, are not elements of $H^{1}(\mathbb{R}^{n})$. Therefore, the variational functional corresponded to (\ref{single}) is not well defined on $H^{1}(\mathbb{R}^{n})$ and the general variational methods are invalid.

On the other hand, since $\lambda>0$ is contained in the essential spectrum of $-\Delta$, the uniqueness of solutions for Helmholtz equation also become complicate. Since different solutions vanish
at infinity, it is necessary to introduce additional conditions at
infinity to determine the unique solution, that is the Sommerfeld (outgoing) radiation condition
\begin{equation*}
u(r)=O(r^{\frac{1-n}{2}}),~~~~\frac{\partial u}{\partial r}+i\sqrt{\lambda}u=o(r^{\frac{1-n}{2}}),
\end{equation*}
see \cite{Rellich1943}. One can also use the ``Limiting Absorption Principle'' to characterize the solution of (\ref{single}). Namely, by constructing the auxiliary problems
\begin{equation}\label{auxiliary}
-\Delta u_{\varepsilon}-(\lambda+i\varepsilon)u_{\varepsilon}=f(x,u_{\varepsilon})~~~\mathrm{in}~~\mathbb{R}^{N},
\end{equation}
one can exclude some solutions $u_{\varepsilon}$ without direct physical significance, which increase without bound as $r\longrightarrow\infty$, and for the remains which are equal to zero at infinity can be used to characterize the solution $u$ of (\ref{single}) as we take $\varepsilon\longrightarrow 0^{+}$, see \cite[Section 7]{Tychonov1967}. Actually,
Kenig, Ruiz and Sogge in \cite{Kenig1987} established the boundeness estimate for the corresponding resolvent
\begin{equation}\label{estimate1}
||u_{\varepsilon}||_{L^{q}(\mathbb{R}^{n})}=||(-\Delta-(\lambda+i\varepsilon))^{-1}f||_{L^{q}(\mathbb{R}^{n})}
\leq||f||_{L^{p}(\mathbb{R}^{n})},
\end{equation}
where $\frac{1}{p}-\frac{1}{q}=\frac{2s}{n}$ and $\mathrm{min}\{|\frac{1}{p}-\frac{1}{2}|,|\frac{1}{q}-\frac{1}{2}|\}>\frac{1}{2n}$. Later, Guti\'{e}rrez  \cite[Theorem 6]{Gutierrez2004} generalized this result to the case of  $\frac{2}{n+1}<\frac{1}{p}-\frac{1}{q}<\frac{2}{n}, \frac{1}{p}>\frac{n+1}{2n}, \frac{1}{q}<\frac{n-1}{2n}$. Taking $\varepsilon\longrightarrow 0^{+}$, one can also obtain the boundedness estimate for the resolvent $\mathcal{R}_{\lambda}=(-\Delta-\lambda)^{-1}$,
see \cite[Theorem 6]{Gutierrez2004} or see \cite{Mandel2019-1,Odeh1961,Radosz2010,Cossetti2021,Mandel2021}. Next, if $f(x,u)=F(x)$ is given integrable function and is independent on $u$, the unique solution then can be directly determined by adding the weaken version of Sommerfeld (outgoing) radiation condition
\begin{equation*}
\mathop{\mathrm{lim}}\limits_{R\longrightarrow+\infty}\frac{1}{R}\int_{B_{R}}|\nabla u-i\sqrt{\lambda}u\widehat{x}|^{2}dx=0,
\end{equation*}
see \cite[Corollary 1]{Gutierrez2004}. And if $f(x,u)$ is dependent on $u$, one can use some fixed point arguments to find the unique solution. Follow this strategy, Guti\'{e}rrez in the celebrated paper \cite{Gutierrez2004}, considered the contour problem for the stationary Ginzburg-Landau equation. By reformulating the Ginzburg-Landau into an integral equation
\begin{equation}\label{interal equation}
u=\varphi+u_{sc}=\widehat{hd\sigma}+((-\Delta+\lambda)^{-1}|u|^{2}u),~~\mathrm{in}~\mathbb{R}^{n},
\end{equation}
the author obtained the nontrivial $L^{q}$ complex solutions via the contraction map argument, where $\varphi=\widehat{hd\sigma}$ is a fourier extension for $h\in L^{2}(S^{n-1},\mathbb{R})$ and it is also called the Herglotz ( incident free) wave function solving the homogeneous Helmholtz equation, while $u_{sc}$ is a scattered wave. We would mention that some smallness condition, $||\varphi||_{L^{4}}(\mathbb{R}^{n})\leq\varepsilon$, is used in the contraction map argument. Later Chen, Ev\'{e}quoz and Weth \cite{Chen2021} established a $L^{\infty}(\mathbb{R}^{n})$ estimate for the resolvent operator, that is
\begin{equation}\label{estimate2}
||u_{\varepsilon}||_{L^{\infty}_{\tau(\alpha)}(\mathbb{R}^{n})}=||(-\Delta-(\lambda+i\varepsilon))^{-1}f||_{L^{\infty}_{\tau(\alpha)}(\mathbb{R}^{n})}
\leq||f||_{L^{\infty}_{\alpha}(\mathbb{R}^{n})},
\end{equation}
where $\tau(\alpha)$ is a parameter with respect to $\alpha$ and $||u||_{L^{\infty}_{\alpha}(\mathbb{R}^{n})}=||\langle x\rangle^{\alpha}u||_{L^{\infty}(\mathbb{R}^{n})}=||(1+|x|^{2})^{\frac{\alpha}{2}}u||_{L^{\infty}(\mathbb{R}^{n})}$.
As a consequence, if $f(x,u)$ is assumed to be the  asymptotically linearly bounded nonlinearities with Lipschitz condition, the authors removed the smallness conditions and obtained the $L^{\infty}(\mathbb{R}^{n})$ solution for a given incident wave $\varphi\in L^{\infty}(\mathbb{R}^{n})$; If $f(x,u)=Q(x)|u|^{p-2}u$ and $Q(x)$ have compact supports with some control of its diameter, the authors obtained a global continuation of the branch $\lambda\longrightarrow u_{\lambda}$ of solutions of the equation
\begin{equation}
u=(-\Delta-\lambda)^{-1}[Q(x)|u|^{p-2}u]+\lambda \varphi.
\end{equation}

For the case $0<s<1$, there is also a growing interest in the fractional Laplacian equation due to various applications involving nonlocal diffusion \cite{Bucur2016}. Since the celebrate work of Caffarelli and Silvestre \cite{Caffarelli2007}, many authors extended the corresponded results of the Laplace equation to the fractional case, see \cite{Chen2020,Cabre2014} and the references therein. We would mention that, especially, Frank, Lenzmann and Silvestre\cite{Frank2015} obtained the uniqueness of radial solutions for the following fractional Laplacian
\begin{equation}
(-\Delta)^{s}+u=|u|^{2^{\ast}_{s}-2}u,~~\mathrm{in}~\mathbb{R}^{n},
\end{equation}
where $2^{\ast}_{s}=\frac{2n}{n-2s}$ is the critical index of fractional Sobolev embedding. While for the fractional Helmholtz equation, there is a few results on it. Recently, Cheng, Li and Yang \cite{Yang2023} proved the equivalence of the homogeneous Helmholtz equation and the fractional Helmholtz equation with arbitrary order, which improves a result of Guan, Murugan and Wei \cite{Wei2023}. This motivates us to consider the fractional Helmholtz equation with linear or superlinear term.

To obtain the complex valued solution for (\ref{main system 1}), we need the analogous boundedness estimate for the resolvent operator $((-\Delta)^{s}-\lambda)^{-1}$. Luckily, based on the work of Kenig, Ruiz and Sogge \cite{Kenig1987}, Huang, Yao, Zheng \cite[Theorem 1.4]{Yao2016} established the boundeness estimate for fractional resolvent, that is
\begin{equation}\label{estimate3}
||u_{\varepsilon}||_{L^{q}(\mathbb{R}^{n})}\leq C_{p,q}||((-\Delta)^{s}-(\lambda+i\varepsilon))^{-1}u_{\varepsilon}||_{L^{q}(\mathbb{R}^{n})},
\end{equation}
where $\frac{n}{n+1}<s<\frac{n}{2}$, $\frac{1}{p}-\frac{1}{q}=\frac{2s}{n}$ and $\mathrm{min}\{|\frac{1}{p}-\frac{1}{2}|,|\frac{1}{q}-\frac{1}{2}|\}>\frac{1}{2n}$. In \cite{Yao2016}, the resolvent of the fractional power of the negative Laplacian is expressed by some sums of the one order resolvent, hence the resolvent estimate has been reduced into the uniform estimate for the Schwartz kernel of the resolvent. However, this method mainly relies on the Stein's oscillatory integral theorem, see \cite[Lemma 2.4]{Kenig1987}, it follows that the exponent in (\ref{estimate3}) maybe not the optimal. By reviving the method of Guti\'{e}rrez \cite{Gutierrez2004}, we obtain the following boundedness estimate for resolvent operator $((-\Delta)^{s}-\lambda)^{-1}$. More precisely, define $\mathcal{R}_{\lambda,\varepsilon}^{s}f=((-\Delta)^{s}-(\lambda+i\varepsilon))^{-1}f$, we have the main theorem as follow.
\begin{thm}\label{thm1}
Let $n\geq 3$, $\frac{n}{n+1}\leq s<\frac{n}{2}$ and  $1<p<q<\infty$ are Lebesgue exponents
satisfying
\begin{equation}
\frac{2}{n+1}\leq \frac{1}{p}-\frac{1}{q}\leq\frac{2s}{n},~\frac{1}{p}>\frac{n+1}{2n},~\frac{1}{q}<\frac{n-1}{2n},
\end{equation}
then there is a uniform constant $C_{p,q}<\infty$ such that for any $\varepsilon>0$
\begin{equation}\label{estimate 1}
||u_{\varepsilon}||_{L^{q}(\mathbb{R}^{n})}=C_{p,q}\lambda^{\frac{n}{2s}(\frac{1}{p}-\frac{1}{q}-1)}||\mathcal{R}_{\lambda,\varepsilon}^{s}f||_{L^{q}(\mathbb{R}^{n})}
\leq||f||_{L^{p}(\mathbb{R}^{n})}.
\end{equation}
On the other hand, if $0<s<\frac{n}{n+1}$, then no such uniform estimates exist. Particularly, there exists a linear operator $\mathcal{R}^{s}_{\lambda}:\mathcal{S}\longrightarrow \mathcal{S}'$ given by
\begin{equation}\label{limit}
\langle \mathcal{R}^{s}_{\lambda}f,g\rangle=\mathop{\mathrm{lim}}\limits_{\varepsilon\longrightarrow 0}\int_{\mathbb{R}^{N}}[\mathcal{R}_{\lambda,\varepsilon}^{s}f](x)g(x)dx,~~\mathrm{for~all}~f,g\in\mathcal{S}.
\end{equation}
\end{thm}

Furthermore, define $D_{x}^{s}u_{\varepsilon}$ by $\widehat{D_{x}^{s}u_{\varepsilon}}(\xi)=|\xi|^{s}\widehat{u_{\varepsilon}}(\xi)$, we then obtain the local $L^{2}(\mathbb{R}^{n})$-estimate for $u_{\varepsilon}$ and $D_{x}^{s}u_{\varepsilon}$.
\begin{thm}\label{local}
Assume that $\frac{n}{n+1}\leq s<\frac{n}{2}$. Let $u_{\varepsilon}=\mathcal{R}_{\lambda,\varepsilon}^{s}f$, then there exists a constant $C$, indepenent of $\lambda$ and $\varepsilon$, such that
\begin{equation}\label{local u}
\mathop{\mathrm{sup}}\limits_{x_{0},R\geq1/\sqrt{\lambda}}\Big(\frac{1}{R}\int_{B(x_{0},R)}|u_{\varepsilon}(x)|^{2}dx  \Big)^{1/2}\leq C\lambda^{\frac{n}{2}(\frac{1}{p}-\frac{1}{2})-\frac{3}{4}}||f||_{L^{p}(\mathbb{R}^{n})},
\end{equation}
whenever $\frac{1}{n+1}\leq\frac{1}{p}-\frac{1}{2}<\frac{s}{2}$ for $n=3,4$ or $\frac{1}{n+1}\leq\frac{1}{p}-\frac{1}{2}<\frac{s}{n}$ for $n\geq5$. Moreover,
\begin{equation}\label{local nablau}
\mathop{\mathrm{sup}}\limits_{x_{0},R\geq1/\sqrt{\lambda}}\Big(\frac{1}{R}\int_{B(x_{0},R)}
|D_{x}^{s} u_{\varepsilon}(x)|^{2}dx  \Big)^{1/2}\leq C\lambda^{\frac{n}{2}(\frac{1}{p}-\frac{1}{2})-\frac{1}{4}}||f||_{L^{p}(\mathbb{R}^{n})},
\end{equation}
whenever $\frac{1}{n+1}\leq\frac{1}{p}-\frac{1}{2}<\frac{s}{n}$ for $n\geq3$.
\end{thm}
Combining these estimates and the limiting absorption principle, we obtain the solutions of $(-\Delta)^{s} u+\lambda u=f$ characterized by the Sommerfeld radiation condition
\begin{equation}
\mathop{\mathrm{lim}}\limits_{R\longrightarrow\infty}\int_{B_{R}}|D_{x}^{s}u-iu\widehat{x}|^{2}dx=0,
\end{equation}
where $\widehat{x}=\frac{x}{|x|}$. Moreover, based on these priori estimates for resolvent, we can obtain the complex solutions for (\ref{main system 1}).
\begin{thm}\label{thm3}
(i) Let $n\geq3$, $\frac{n}{n+1}\leq s<\frac{n+1}{4}$ and $ns-n-s>0$. Moreover, let
\begin{equation}
\left \{
\begin{aligned}
&\frac{n(t-1)}{2s}<q<\frac{2nt}{n+1},~~&&\mathrm{if}~~\frac{(n+1)^{2}}{(n-1)^{2}}<t<\frac{n+1}{n+1-4s},\\
&\frac{n(t-1)}{2s}<q<\frac{(n+1)(t-1)}{2},~~&&\mathrm{if}~~\frac{(n-1+4s)}{(n-1)}<t<\frac{(n+1)^{2}}{(n-1)^{2}},\\
&\frac{2n}{n-1}<q<\frac{(n+1)(t-1)}{2},~~&&\mathrm{if}~~\frac{n^{2}+4n-1}{n^{2}-1}<t<\frac{n-1+4s}{n-1};\\
\end{aligned}
\right.
\end{equation}
(ii) Let $n\geq3$, $\frac{n}{n+1}\leq s<\frac{n+1}{4}$ and $ns-n-s<0$. Moreover, let
\begin{equation}
\left \{
\begin{aligned}
&\frac{n(t-1)}{2s}<q<\frac{2nt}{n+1},~~&&\mathrm{if}~~\frac{n-1+4s}{n-1}<t<\frac{n+1}{n+1-4s},\\
&\frac{2n}{n-1}<q<\frac{2nt}{n+1},~~&&\mathrm{if}~~\frac{(n+1)^{2}}{(n-1)^{2}}<t<\frac{n-1+4s}{n-1},\\
&\frac{2n}{n-1}<q<\frac{(n+1)(t-1)}{2},~~&&\mathrm{if}~~\frac{n^{2}+4n-1}{n^{2}-1}<t<\frac{(n+1)^{2}}{(n-1)^{2}};\\
\end{aligned}
\right.
\end{equation}
(iii) Let $n=3,4$, $\frac{n+1}{4}<s<\frac{2n^{2}}{(n+1)^{2}}$. Moreover, let
\begin{equation}
\left \{
\begin{aligned}
&\frac{n(t-1)}{2s}<q<\frac{2nt}{n+1},~~&&\mathrm{if}~~\frac{(n+1)^{2}}{(n-1)^{2}}<t<+\infty,\\
&\frac{n(t-1)}{2s}<q<\frac{(n+1)}{t-1},~~&&\mathrm{if}~~\frac{n}{n-2s}<t<\frac{(n+1)^{2}}{(n-1)^{2}},\\
&t<q<\frac{(n+1)(t-1)}{2},~~&&\mathrm{if}~~\frac{2n}{n-1}<t<\frac{n}{n-2s},\\
&\frac{2n}{n-1}<q<\frac{(n+1)(t-1)}{2},~~&&\mathrm{if}~~\frac{n^{2}+4n-1}{n^{2}-1}<t<\frac{2n}{n-1};\\
\end{aligned}
\right.
\end{equation}
(iv) Let $n=3,4$, $\frac{2n^{2}}{(n+1)^{2}}\leq s<\frac{n}{2}$. Moreover, let
\begin{equation}
\left \{
\begin{aligned}
&\frac{n(t-1)}{2s}<q<\frac{2nt}{n+1},~~&&\mathrm{if}~~\frac{n}{n-2s}<t<+\infty,\\
&t<q<\frac{2nt}{n+1},~~&&\mathrm{if}~~\frac{(n+1)^{2}}{(n-1)^{2}}<t<\frac{n}{n-2s},\\
&t<q<\frac{(n+1)(t-1)}{2},~~&&\mathrm{if}~~\frac{2n}{n-1}<t<\frac{(n+1)^{2}}{(n-1)^{2}},\\
&\frac{2n}{n-1}<q<\frac{(n+1)(t-1)}{2},~~&&\mathrm{if}~~\frac{n^{2}+4n-1}{n^{2}-1}<t<\frac{2n}{n-1};\\
\end{aligned}
\right.
\end{equation}
(v) Let $n\geq5$, $\frac{n+1}{4}<s<\frac{n}{2}$. Moreover, let
\begin{equation}
\left \{
\begin{aligned}
&\frac{n(t-1)}{2s}<q<\frac{2nt}{n+1},~~&&\mathrm{if}~~\frac{n}{n-2s}<t<+\infty,\\
&t<q<\frac{2nt}{n+1},~~&&\mathrm{if}~~\frac{2n}{n-1}<t<\frac{n}{n-2s},\\
&\frac{2n}{n-1}<q<\frac{2nt}{n+1},~~&&\mathrm{if}~~\frac{(n+1)^{2}}{(n-1)^{2}}<t<\frac{2n}{n-1},\\
&\frac{2n}{n-1}<q<\frac{(n+1)(t-1)}{2},~~&&\mathrm{if}~~\frac{n^{2}+4n-1}{n^{2}-1}<t<\frac{(n+1)^{2}}{(n-1)^{2}}.\\
\end{aligned}
\right.
\end{equation}
Assume that $f(x,u)=|u|^{t-1}u$, then for any given $\varphi\in L^{q}(\mathbb{R}^{n})$ of the homogeneous Helmholtz equation $(-\Delta)^{s}\varphi-\lambda\varphi=0$ with $||\varphi||_{L^{q}(\mathbb{R}^{n})}\leq\varepsilon$, there exists $a=a(||\varphi||_{L^{q}(\mathbb{R}^{n})})$ such that (\ref{main system 1}) has a unique solution $u=\mathcal{R}^{s}_{\lambda}(|u|^{t-1}u)+\varphi\in L^{q}(\mathbb{R}^{n})$ satisfying
\begin{equation*}
u\in B_{a}(L^{q}(\mathbb{R}^{n}))=\{u:\mathbb{R}^{n}\longrightarrow \mathbb{C}|~||u||_{L^{q}(\mathbb{R}^{n})}\leq a\}.
\end{equation*}
\end{thm}
\begin{Rem}
Indeed, the fractional Helmholtz equation are very closed to the fractional Ginzburg-Landau equations and fractional Allen-Cahn equtaion, see \cite{Millot2015,Wang2016} and the references therein. Therefore, by the similar analysis as in Guti\'{e}rrez \cite{Gutierrez2004} and the boundedness estimate in Ma \cite{Ma2016}, some existence result can also be obtained for the fractional Ginzburg-Landau equations.
\end{Rem}
To remove the smallness condition $||\varphi||_{L^{q}(\mathbb{R}^{n})}\leq\varepsilon$,  we also establish the following boundedness estimate for resolvent operator $((-\Delta)^{s}-\lambda)^{-1}$.
\begin{thm}\label{thm5}
Let $n\geq3$, $1\leq s\leq\frac{n}{2}$, $\alpha>\frac{n+1}{2}$, and $\tau(\alpha)$ be defined by
\begin{equation}
\tau(\alpha)=\left \{
\begin{aligned}
&\alpha-\frac{n+1}{2}, ~&&\mathrm{if}~\frac{n+1}{2}<\alpha<n,\\
&\frac{n+1}{2}, ~ &&\mathrm{if}~\alpha\geq n.\\
\end{aligned}
\right.
\end{equation}
Then we have
\begin{equation}\label{estimate4}
\kappa_{\alpha}:=\mathrm{sup}\{||\mathcal{R}^{s}_{\lambda}f||_{L^{\infty}_{\tau(\alpha)}(\mathbb{R}^{n})}:f\in L^{\infty}_{\alpha}(\mathbb{R}^{n}),||f||_{L^{\infty}_{\alpha}(\mathbb{R}^{n})}=1\}<\infty.
\end{equation}
So $\mathcal{R}^{s}_{\lambda}$ defines a bounded linear map $L^{\infty}_{\alpha}(\mathbb{R}^{n})\longrightarrow L^{\infty}_{\tau(\alpha)}(\mathbb{R}^{n})$. Moreover, the resolvent operator defines a compact linear map $\mathcal{R}^{s}_{\lambda}:L^{\infty}_{\alpha}(\mathbb{R}^{n})\longrightarrow L^{\infty}(\mathbb{R}^{n})$.
\end{thm}

Correspondingly, we have the following existence results for (\ref{main system 1}).
\begin{thm}\label{thm6}
Let $n\geq3$, $1\leq s\leq\frac{n}{2}$. For some $\alpha>\frac{n+1}{2}$, let $f:\mathbb{R}^{n}\times\mathbb{C}\longrightarrow\mathbb{C}$ be a
continuous function satisfying
\begin{equation}\label{condition1}
\mathop{\mathrm{sup}}\limits_{|u|\leq M,x\in\mathbb{R}^{n}}\langle x\rangle^{\alpha}|f(x,u)|<\infty~~\mathrm{for~all}~M>0.
\end{equation}
Moreover, suppose that the nonlinearity is of the form $f(x,u)\leq Q(x)|u|+b(x)$ with $Q,b\in L^{\infty}_{\alpha}(\mathbb{R}^{n},\mathbb{R})$ and $||Q||_{L_{\alpha}^{\infty}(\mathbb{R}^{n})}$. Then, for any given solution $\varphi\in L^{\infty}(\mathbb{R}^{n})$ of the homogeneous Helmholtz equation $(-\Delta)^{s}\varphi+\lambda \varphi=0$, (\ref{main system 1}) admits a solution $u=\mathcal{R}^{s}_{\lambda}(f(x,u))+\varphi\in L^{\infty}(\mathbb{R}^{n})$. Particulary, if $f(x,u)$ satisfies the Lipschitz condition
\begin{equation}\label{lips}
l_{\alpha}:=\mathrm{sup}\Big\{\langle x\rangle^{\alpha}\Big|\frac{f(x,u)-f(x,v)}{u-v}\Big|:u,v\in\mathbb{R},x\in\mathbb{R}^{n}\Big\}\leq \frac{1}{\kappa_{\alpha}},
\end{equation}
then the solution is unique.
\end{thm}

\begin{Rem}
The case ($f_{1}$) in \cite{Chen2021} 
is more complicate, someone need more nonexistence result for fractional Helmholtz equation to deal with it.
\end{Rem}

As for equation (\ref{main system 1}) with superlinear term, we can obtain the similar result as Theorem \ref{thm3}.
\begin{thm}\label{thm8}
Let $n\geq3$, $1\leq s\leq\frac{n}{2}$. For some $\alpha>\frac{n+1}{2}$, let $f(x,u):\mathbb{R}^{n}\times\mathbb{C}\longrightarrow\mathbb{C}$ be a continuous function satisfies (\ref{condition1}). Suppose that $f(x,\cdot)$ is real defferentiable for every $x\in\mathbb{R}^{n}$ and $f':=\partial_{u}f:\mathbb{R}^{n}\times\mathbb{C}\longrightarrow\mathcal{L}_{\mathbb{R}}(\mathbb{C},\mathbb{C})$ is a continuous function satisfying
\begin{equation}
\mathop{\mathrm{sup}}\limits_{|u|\leq M,x\in\mathbb{R}^{n}}\langle x\rangle ||f'(x,u)||_{\mathcal{L}_{\mathbb{R}}(\mathbb{C},\mathbb{C})}<\infty.
\end{equation}
Moreover, suppose that $f(x,0)=0$ and $f'(x,0)=0\in\mathcal{L}_{\mathbb{R}}(\mathbb{C},\mathbb{C})$ for all $x\in\mathbb{R}^{n}$. Then there exists open neighborhoods $U,V\subset L^{\infty}(\mathbb{R}^{n})$ of zero with the property that for every $\varphi\in V$ there exists a unique solution $u=u_{\varphi}\in U$ of (\ref{main system 1}). Moreover, the map $V\longrightarrow U$, $u\longrightarrow u_{\varphi}$ is of class $C^{1}$.
\end{thm}

Set $f(x,u)=Q(x)|u|^{p-2}u$ with $p>2$ and $Q\in L^{\infty}_{\alpha}(\mathbb{R}^{n})$ for some $\alpha>\frac{n+1}{2}$, then $f(x,u)$ is a special example that satisfying the conditions of Theorem \ref{thm8}. Therefore, for given $\varphi\in L^{\infty}(\mathbb{R}^{n})$, there exists $\epsilon>0$ and a unique local branch $(\epsilon,\epsilon)\longrightarrow L^{\infty}(\mathbb{R}^{n}),\lambda\longrightarrow u_{\lambda}$ of solutions the equation
\begin{equation}\label{branch}
u=\mathcal{R}_{\lambda}^{s}(Q|u|^{p-2}u)+\lambda\varphi~~~\mathrm{in}~L^{\infty}(\mathbb{R}^{n}).
\end{equation}

\begin{Rem}
To establish the existence of a global continuation of this local branch, someone also need the nonexistence for fractional Helmholtz equation, see the classical case in \cite[Section 4]{Chen2021}.
Even though one may assume that the stronger condition, that is  $Q\in L^{\infty}_{c}(\mathbb{R}^{n},\mathbb{R})\setminus\{0\}$ with some control on its diameters, the specific form of the Green function for fractional Helmholtz operator is not clear.
\end{Rem}

We would mention that some real valued solutions $u=\mathrm{Re}(\mathcal{R}_{\lambda}f)$ for (\ref{single}) haven also been detected by many authors. Since the real-valued solutions is only a real part of convolution integral, it easily follows that $u=0$ is an isolated solution of (\ref{single}) in $L^{p}(\mathbb{R}^{n})$, and thus the nontrivial solutions cannot be found by a contraction mapping argument. Therefore, based on the boundedness estimate of the resolvent operator $\mathcal{R}_{\lambda}$ with $\frac{2(n+1)}{n-1}\leq p\leq \frac{2n}{n-2}$, Ev\'{e}quoz and Weth \cite{Evequoz2015} (see \cite{Evequoz2016} for n=2) set up a dual variational framework for (\ref{single}). Correspondingly, the nontrivial real-valued solutions of equation (\ref{single}) with $f(x,u)=Q(x)|u|^{p-2}u$ are detected via the mountain pass argument, where $Q(x)$ is a periodic or decay weight function,  see also  \cite{Evequoz2017-2,Evequoz2015-1,Mandel2021-1,Evequoz2017-1,Evequoz2020,Mandel2019} for the other cases. Specially, Ev\'{e}quoz and Weth \cite{Evequoz2017-3} obtained the positive solution for (\ref{branch}). And by setting $h(\xi)=\overline{h(-\xi)}$ in (\ref{interal equation}), Mandel \cite{Mandel2019} revived Guti\'{e}rrez' fixed point approach and detected the continua of small real-valued solutions
of (\ref{single}) for a larger class of nonlinearities. In this paper, we also consider the real valued solutions for (\ref{main system 1}).
\begin{thm}\label{thm10}
Let $n\geq 3$, $\frac{n}{n+1}<s<\frac{n}{2}$, $\frac{2(n+1)}{n-1}<p<\frac{2n}{n-2s}$, and let $Q\in L^{\infty}(\mathbb{R}^{n})$, $Q\geq0$, $Q\not\equiv0$ satisfy $\mathop{\mathrm{lim}}\limits_{|x|\longrightarrow\infty}Q(x)=0$. Then problem $u=\mathrm{Re}(\mathcal{R}_{\lambda}^{s}(Q(x)|u|^{p-2}u))$ admits a sequence of pairs $\pm u_{n}$ of solutions such that $u_{n}\in L^{p}(\mathbb{R}^{n})$ with $||u_{n}||_{L^{p}(\mathbb{R}^{n})}\longrightarrow \infty$ as $n\longrightarrow \infty$. Particularly, $u_{n}\in W^{2s,q}(\mathbb{R}^{n})\cap \mathcal{C}^{1,\alpha}(\mathbb{R}^{n})$ for $q\in[p,\infty)$ and $\alpha\in(0,1)$.
\end{thm}

\begin{thm}\label{thm11}
Let $n\geq 3$, $\frac{n}{n+1}<s<\frac{n}{2}$, $\frac{2(n+1)}{n-1}<p<\frac{2n}{n-2s}$, and let $Q\in L^{\infty}(\mathbb{R}^{n})$, $Q\geq0$, $Q\not\equiv0$ be $\mathbb{Z}^{n}$-periodic. Then problem  $u=\mathrm{Re}(\mathcal{R}_{\lambda}^{s}(Q(x)|u|^{p-2}u))$ admits a nontrivial solution such that $u\in L^{p}(\mathbb{R}^{n})$. Particularly, $u_{n}\in W^{2s,q}(\mathbb{R}^{n})\cap \mathcal{C}^{1,\alpha}(\mathbb{R}^{n})$ for $q\in[p,\infty)$ and $\alpha\in(0,1)$.
\end{thm}

\begin{Rem}
The far filed estimate for the real valued solutions need more information of the specific form on the asymptotic expansions for the Green function of fractional Helmholtz operator.
\end{Rem}

Let us now briefly explain our approach and the organization of the paper. In Sections \ref{pre}, we first derive the limiting absorption principle for fractional Helmholtz operator, i.e., the $L^{p}(\mathbb{R}^{n})$ and $L^{\infty}(\mathbb{R}^{n})$ estimate for resolvent. In Section \ref{complex}, we prove the existence of the complex valued solutions for (\ref{main system 1}). In Section \ref{non}, we derive a nonvanishing property related to the resolvent which is a key ingredient in the proof of Theorem 10. In Section \ref{var}, we lift the regularity of the solutions for a priori equation, and then we obtain some compactness for the resolvent operator. With the help of this property, we set up a dual variational framework for the problem $u=\mathrm{Re}(\mathcal{R}_{\lambda}^{s}f)$. In section \ref{real}, we obtain the existence of real valued solutions for (\ref{main system 1}).

%
\setlength{\parindent}{2em}
\section{Limiting absorption principle for fractional Helmholtz operator}\label{pre}
Since the resolvent estimate for fractional Hemholtz operator is almost similar to the classical Helmholtz operator, hence in the first subsection, we recall and compare the delicate difference between the resolvent estimate in Guti\'{e}rrez \cite{Gutierrez2004} and the integral estimate in Chen, Ev\'{e}quoz and Weth \cite{Chen2021}.
\subsection{Green function for Helmholtz operator}
 Let $\lambda,\varepsilon>0$. Then the operator $-\Delta-(\lambda+i\varepsilon):H^{2}(\mathbb{R}^{n})\subset L^{2}(\mathbb{R}^{n})\longrightarrow L^{2}(\mathbb{R}^{n})$ is an isomorphism. Moreover, for any $f$ from the Schwartz space $\mathcal{S}$ its inverse is given by
\begin{equation*}
\mathcal{R}_{\lambda,\varepsilon}f(x):=[-\Delta-(\lambda+i\varepsilon)]^{-1}f(x)=(2\pi)^{-\frac{N}{2}}\int_{\mathbb{R}^{n}}e^{ix\cdot\xi}\frac{\widehat{f}(\xi)}{|\xi|^{2}-(\lambda+i\varepsilon)}d\xi.
\end{equation*}
According to the Limiting Absorption Principle of Guti\'{e}rrez \cite{Gutierrez2004} (see also  \cite{Gelfand1964}) that there exists a linear operator $\mathcal{R}_{\lambda}:\mathcal{S}\longrightarrow\mathcal{S}'$ given by
\begin{equation*}
\langle\mathcal{R}_{\lambda}f,g\rangle:=\mathop{\mathrm{lim}}\limits_{\varepsilon\longrightarrow0}\int_{\mathbb{R}^{n}}
[\mathcal{R}_{\lambda,\varepsilon}f](x)g(x)dx=\int_{\mathbb{R}^{n}}[\Phi_{\lambda}\ast f](x)g(x)dx~~\mathrm{for}~f,g\in\mathcal{S}
\end{equation*}
with
\begin{equation*}
\Phi_{\lambda}(x):=(2\pi)^{-\frac{n}{2}}\mathcal{F}^{-1}((|\xi|^{2}-\lambda-i0)^{-1})(x)
=\lambda^{\frac{n-1}{2}}\Phi_{1}(\sqrt{\lambda}x)=\frac{i}{4}(\frac{\lambda}{4\pi^{2}|x|^{2}})^{\frac{2-n}{4}}H_{\frac{n-2}{2}}^{(1)}(\sqrt{\lambda}|x|)
\end{equation*}
for $x\in\mathbb{R}^{n}\setminus\{0\}$, where $H^{(1)}_{\frac{n-2}{2}}$ is the Hankel function of the first kind of order $\frac{n-2}{2}$. Here we use the notation form \cite{Gelfand1964}, which also allows us briefly write
\begin{equation*}
\mathcal{R}_{\lambda}f:=\mathcal{F}^{-1}\left((|\xi|^{2}-\lambda-i0)^{-1}\widehat{f}\right)~~\mathrm{for}~f\in\mathcal{S}.
\end{equation*}
For $H^{(1)}_{\frac{n-2}{2}}$ we have the asymptotic expansions
\begin{equation*}
H^{(1)}_{\frac{n-2}{2}}(s)=\left \{
\begin{aligned}
&\sqrt{\frac{2}{\pi s}}e^{i(s-\frac{n-1}{4}\pi)}[1+O(s^{-1})], ~&&\mathrm{as}~~s\longrightarrow\infty,\\
&-\frac{i\Gamma(\frac{n-2}{2})}{\pi}\Big(\frac{2}{s}\Big)^{\frac{N-2}{2}}[1+O(s)], ~&&\mathrm{as}~~s\longrightarrow 0^{+},
\end{aligned}
\right.
\end{equation*}
(see e.g. \cite[Formulas (5.16.3)]{Lebedev1972}), so there exists a constant $C_{0}>0$ such that
\begin{equation}\label{Phi}
|\Phi_{\lambda}(x)|\leq \left \{
\begin{aligned}
&C_{0}\mathrm{max}\{|x|^{2-n},|x|^{\frac{1-n}{2}}\}~~\mathrm{for}~~x\in\mathbb{R}^{n}\setminus\{0\}, ~&&N\geq3,\\
&C_{0}\mathrm{min}\{1+|\mathrm{log}~|x||,|x|^{-\frac{1}{2}}\}~~\mathrm{for}~~x\in\mathbb{R}^{n}\setminus\{0\}, ~&&N=2.
\end{aligned}
\right.
\end{equation}

As we can see, $\Phi_{\lambda}$ is a Green function of Helmholtz operator but without the uniform bounded estimate. Actually, let $z\in\mathbb{C}$, and let $\Phi_{z}(x)$ be the Fourier transform of the multipliers $m_{\lambda}=\frac{1}{|\xi|^{2}+z}$, then, as $arg~z\notin [-\frac{\pi}{2},\frac{\pi}{2}]$, some exponential decay properties of $H_{\frac{N-2}{2}}^{1}(\sqrt{|\xi|^{2}z})$ is not uniform, see more details in Kenig, Ruiz, Sogge \cite{Kenig1987}. Based on these essential properties of the Green function, an oscillatory integral theorem of Stein \cite{Stein1986} has  been used to deal with this problem, see \cite{Kenig1987}. This method also be used in the case of fractional Helmholtz operator, see Huang, Yao and Zheng \cite{Yao2016}. However, these result may be not the optimal. Hence, a cut-off skill and harmonic analysis method has been proposed by Guti\'{e}rrez \cite{Gutierrez2004}. As a consequence, a more precisely version of the resolvent estimate is established. This method also be used to establish the nonvanishing lemma in Ev\'{e}quoz and Weth \cite{Evequoz2015}. While, in the paper of Chen, Ev\'{e}quoz and Weth, they assumed that the nonlinearities $f(x)$ belong to a stronger integrable space $L^{\infty}_{\alpha}(\mathbb{R}^{n})$, that is some functions satisfying decay condition, and then they obtained the bounded estimate for the resolvent with the help of the weight term $\langle x\rangle^{\alpha}=(1+|x|^{2})^{\frac{\alpha}{2}}$.

\subsection{$L^{p}(\mathbb{R}^{n})$ estimate for resolvent of fractional Helmholtz operator}
Follow the idea of Guti\'{e}rrez \cite{Gutierrez2004}, we give the proof of Theorem \ref{thm1}.

\begin{proof}[\bf Proof of Theorem \ref{thm1}]
For any $f(x)\in \mathcal{S}$, the Schwartz space, by using the Fourier transform, the solution
$u_{\varepsilon}=\mathcal{R}_{\lambda,\varepsilon}^{s}f$
then can be written as (here we drop the subscript $\varepsilon$ from notation)
\begin{equation}\label{fourier}
u(x)=c\int_{\mathbb{R}^{n}}e^{ix\cdot\xi}\frac{1}{|\xi|^{2s}-(\lambda+i\varepsilon)}\widehat{f(\xi)}d\xi.
\end{equation}
Since the problem is invariant under dilation and rotation we restrict ourselves to the case $\lambda=1$ and
\begin{equation}
\mathrm{supp}\widehat{f}\subseteq\{\xi=(\overline{\xi},\xi_{n})\in\mathbb{R}^{n-1}\times\mathbb{R}:|\overline{\xi}|<\xi_{n}/6,\xi_{n}>0\}.
\end{equation}
Take a radial cut-off function $\phi\in C_{c}^{\infty}(\mathbb{R})$, with $\mathrm{supp}\phi\subseteq[0,3/4]$, $\phi(t)=1$ if $t\in[0,5/8]$, and $0\leq\phi\leq1$. Then, we can split the multiplier $m(\xi)=(|\xi|^{2s}-(1+i\varepsilon))^{-1}$ in (\ref{fourier}) into multiplier $m_{i}$, $i=1,2,3$, in the following way
\begin{equation}
m_{1}(\xi)=\phi(|\xi|)m(\xi),~~m_{2}(\xi)=(1-\phi(|\xi|/2))m(\xi),~~m_{3}=m-(m_{1}+m_{2}).
\end{equation}
Let $u_{i}$ such that $\widehat{u_{i}}(\xi)=m_{i}(\xi)\widehat{f}(\xi)$ or, equivalently, $u_{i}=M_{i}\ast f$ where $M_{i}$ denotes the Fourier transform of $m_{i}$ (i.e. $\widehat{M_{i}}=m_{i}$), it suffices to prove inequality (\ref{estimate 1}) for each $u_{i}$.

Firstly, it is easy to check that function $u_{i}$, $i=1,2$, are pointwise majorized by the Bessel
potential $J^{2s}f=(I-\Delta)^{s}f$. Hence the Fractional Integral Theorem yields
\begin{equation}\label{young}
||u_{i}||_{L^{q}(\mathbb{R}^{n})}=||M_{i}\ast f||_{L^{q}(\mathbb{R}^{n})}\leq C||J^{2s}f||_{L^{q}(\mathbb{R}^{n})}\leq C||f||_{L^{p}(\mathbb{R}^{n})},
\end{equation}
when $0\leq\frac{1}{p}-\frac{1}{q}\leq\frac{2s}{n}$, $q\neq\infty$ and $p\neq1$.

It remains $u_{3}$ to be estimated. Without loss of generality, we assume that $\mathrm{supp}\widehat{f}$ is contained in a neighbourhood of the support of $m_{3}$, i.e. $||\xi|-1|<1/2$. We then claim that
\begin{equation}
|M_{3}(x)|\leq\frac{C}{(1+|x|)^{\frac{n-1}{2}}},~~a.e.~x\in\mathbb{R}^{n},~\varepsilon>0.
\end{equation}
Indeed, from the definition of $M_{3}$ (i.e. $\widehat{M_{3}}=m_{3}$) we can write
\begin{equation}
M_{3}(x)=M_{3}(\overline{x},x_{n})=\int_{\mathbb{R}^{N-1}}e^{i\overline{x}\cdot\overline{
\xi}}\Big(\int_{\mathbb{R}}e^{ix_{n}\xi_{n}}\frac{\psi(\xi)}{|\xi|^{2s}-(1+i\varepsilon)}d\xi_{n}   \Big)d\overline{\xi},
\end{equation}
where $\psi(\xi)=\phi(|\xi|/2)-\phi(|\xi|)$ is a compactly supported smooth function and
\begin{equation}
\mathrm{supp}\psi\subset \{(\overline{\xi},\xi_{n}):|\overline{\xi}|<1/4,3/8<\xi_{n}<3/2\}.
\end{equation}
Using the change of variables $\eta_{n}=\xi_{n}-(1-|\overline{\xi}|^{2})^{1/2}$, i.e. $\xi_{n}=\eta_{n}+(1-|\overline{\xi}|^{2})^{1/2}$ in the above identity, the kernel $M_{3}$ can be rewritten as
\begin{equation}
\begin{aligned}
&M_{3}(x)=M_{3}(\overline{x},x_{n})=\int_{\mathbb{R}^{n-1}}e^{i\overline{x}\cdot\overline{
\xi}}\Big(\int_{\mathbb{R}}e^{ix_{n}\xi_{n}}\frac{\psi(\xi)}{|\xi|^{2s}-(1+i\varepsilon)}d\xi_{n}   \Big)d\overline{\xi}\\
&=\int_{\mathbb{R}^{n-1}}e^{i\overline{x}\cdot\overline{
\xi}}\Big(\int_{\mathbb{R}}e^{ix_{n}\xi_{n}}\frac{\psi(\overline{\xi},\xi_{n})}{|\xi|^{2s}-(1+i\varepsilon)}d\xi_{n}   \Big)d\overline{\xi}\\
&=\int_{\mathbb{R}^{n-1}}e^{i\overline{x}\cdot\overline{
\xi}}\Big(\int_{\mathbb{R}}e^{ix_{n}(\eta_{n}+(1-|\overline{\xi}|^{2})^{1/2})}\frac{\psi(\overline{\xi},\eta_{n}+(1-|\overline{\xi}|^{2})^{1/2})}{|\xi|^{2s}-(1+i\varepsilon)}d(\eta_{n}+(1-|\overline{\xi}|^{2})^{1/2})   \Big)d\overline{\xi}\\
&=\int_{\mathbb{R}^{n-1}}e^{i\overline{x}\cdot\overline{
\xi}+ix_{n}(1-|\overline{\xi}|^{2})^{1/2}}\Big(\int_{\mathbb{R}}e^{ix_{n}\eta_{n}}\frac{\psi(\overline{\xi},\eta_{n}+(1-|\overline{\xi}|^{2})^{1/2})}
{[|\overline{\xi}|^{2}+(\eta_{n}+(1-|\overline{\xi}|^{2})^{1/2})^{2}]^{s}-(1+i\varepsilon)}d(\eta_{n}+(1-|\overline{\xi}|^{2})^{1/2})   \Big)d\overline{\xi}\\
&=\int_{\mathbb{R}^{n-1}}e^{i\overline{x}\cdot\overline{
\xi}+ix_{n}(1-|\overline{\xi}|^{2})^{1/2}}\Big(\int_{\mathbb{R}}e^{ix_{n}\eta_{n}}\frac{\psi(\overline{\xi},\eta_{n}+(1-|\overline{\xi}|^{2})^{1/2})}
{[\eta_{n}^{2}+2\eta_{n}(1-|\overline{\xi}|^{2})^{1/2}+1]^{s}-(1+i\varepsilon)}d(\eta_{n}+(1-|\overline{\xi}|^{2})^{1/2})   \Big)d\overline{\xi}\\
&=\int_{\mathbb{R}^{n-1}}e^{i\overline{x}\cdot\overline{
\xi}+ix_{n}(1-|\overline{\xi}|^{2})^{1/2}}\Big(\int_{\eta_{n}}e^{ix_{n}\eta_{n}}
\frac{\widetilde{\psi}_{\varepsilon}(\overline{\xi},\eta_{n})}{\eta_{n}^{s}}d\eta_{n}  \Big)d\overline{\xi}\\
&=\int_{\mathbb{R}^{n-1}}e^{i\overline{x}\cdot\overline{
\xi}+ix_{n}(1-|\overline{\xi}|^{2})^{1/2}}\gamma(\overline{\xi},x_{n})d\overline{\xi},\\
\end{aligned}
\end{equation}
where
\begin{equation}
\gamma(\overline{\xi},x_{n})=\int_{\eta_{n}}e^{ix_{n}\eta_{n}}
\frac{\widetilde{\psi}_{\varepsilon}(\overline{\xi},\eta_{n})}{\eta_{n}^{s}}d\eta_{n}
\end{equation}
and $\widetilde{\psi}_{\varepsilon}(\overline{\xi},\eta_{n})$ is a Schwartz function with uniform estimates in $\overline{\xi}$ and $\varepsilon$, $\gamma(\overline{\xi},x_{n})\in\mathcal{C}_{c}^{\infty}(\mathbb{R}^{N-1})$ with support and uniform estimates in the set of $\{\xi:|\overline{\xi}|<\frac{1}{4}\}$.

Denote the phase function of $M_{3}$ by $\sigma(\xi)=\overline{x}\cdot\overline{\xi}+x_{n}(1-|\overline{\xi}|^{2})^{1/2}$, then it is easy to check that $\sigma(\xi)$ has no critical points in the support of $\gamma(\overline{\xi},x_{n})$ when $|\overline{x}|\leq |x_{n}|/3$. Therefore, it suffices consider points $x=(\overline{x},x_{n})$ with $|\overline{x}|\leq|x_{n}|/3$. Notice that $1-|\overline{\xi}|^{2}>c>0$ if $\overline{\xi}\in\mathrm{supp}\gamma(\cdot,t)$, and so
\begin{equation}
\mathrm{det}(\partial_{ij}^{2}\sigma)\geq|x_{n}|^{n-1},~~\forall\xi\in\mathrm{supp}\gamma(\cdot,t).
\end{equation}
Hence, by the stationary phase lemma, see \cite[Theorem 7.7.17]{Hormander1983}, we are lead to
\begin{equation}
|M_{3}(x)|\leq C(1+|x_{n}|)^{-(n-1)/2},~~\forall~\varepsilon>0,
\end{equation}
and consequently, our claims hold. Follow the similar proof in \cite[Theorem 6]{Gutierrez2004}, for $\frac{1}{p}-\frac{1}{q}\geq\frac{1}{n+1}$, $\frac{1}{q}<\frac{n-1}{2n}$, and $\frac{1}{p}>\frac{n+1}{2n}$, we obtain that
\begin{equation}\label{estimate 2}
||u_{3}||_{L^{q}(\mathbb{R}^{n})}=||M_{3}\ast f||_{L^{q}(\mathbb{R}^{n})}\leq C||f||_{L^{p}(\mathbb{R}^{N})}.
\end{equation}
Together with the estimate in (\ref{young}), we obtain the estimate (\ref{estimate 1}).

On the other hand, since the above estimates for $M_{1},M_{2}$ and $M_{3}$ are uniform with respect to $\varepsilon$, hence we have the limiting absorption principle for $\mathcal{R}^{s}_{\lambda,\varepsilon}$, that is (\ref{limit}).
\end{proof}

The proof of Theorem \ref{local} follows the same method in \cite[Theorem 8]{Gutierrez2004}.
\begin{proof}[\bf Proof of Theorem \ref{local}]
Firstly, we prove inequality (\ref{local nablau}). Using the Fourier transform, we can write
\begin{equation}
D_{x}^{s}u(x)=c\int_{\mathbb{R}^{n}}e^{ix\cdot\xi}\frac{\xi^{s}}{|\xi|^{2s}-(\lambda+i\varepsilon)}\widehat{f}(\xi)d\xi.
\end{equation}
Since the problem is invariant under dilations, rotations and translations we restrict ourselves to the case $x_{0}=0,\tau=1$ and
\begin{equation}
\mathrm{supp}\widehat{f}\subset\{\xi=(\overline{\xi,\xi_{n}})\in\mathbb{R}^{n-1}\times\mathbb{R}:|\overline{\xi}|<\xi_{n}/6,\xi_{n}>0\}.
\end{equation}
Let $\phi$ be a radial cut off function such that $\phi\in \mathcal{C}_{c}^{\infty}(\mathbb{R})$ with $\mathrm{supp}\phi\subset[0,3/4]$, $\phi(t)=1$ if $t\in[0,5/8]$, and $0\leq\phi\leq1$. Define $f_{i}$, $i=1,2,3$ by
\begin{equation}
\widehat{f}_{1}(\xi)=\phi(|\xi|)\widehat{f}(\xi),~~\widehat{f}_{2}(\xi)=(1-\phi(|\xi|/2))\widehat{f}(\xi),~~
\widehat{f}_{3}=\widehat{f}-(\widehat{f}_{1}+\widehat{f}_{2}).
\end{equation}
We denote by $u_{i}$ the solution of the equation corresponding to the inhomogeneous term $F_{i}$. Then, it suffices to prove inequality (\ref{local nablau}) for each $D_{x}^{s}u_{i}$.

It is easy to see that $D_{x}^{s}u_{1}$ and $D_{x}^{s}u_{2}$ are pointwise bounded by $J^{s/2}f_{1}$ and $J^{s/2}f_{2}$. Hence, by H\"{o}lder's inequality, Minkowski Integral Inequality, and the known estimates for the Bessel potentials, we are lead to
\begin{equation}
\begin{aligned}
\Big(\frac{1}{R}\int_{B_{R}}|D_{x}^{s}u_{i}|^{2}dx\Big)^{1/2}
&=\Big(\frac{1}{R}\int_{\mathbb{R}^{n}}\chi_{B_{R}}|D_{x}^{s}u_{i}|^{2}dx\Big)^{1/2}\\
&\leq\Big(\frac{1}{R}|B_{R}|^{1-\frac{1}{q}}|||D_{x}^{s}u_{i}|^{2}||_{L^{q}(\mathbb{R}^{n})} \Big)^{1/2}\\
&\leq\Big(\frac{1}{R}|B_{R}|^{1-\frac{1}{q}}||D_{x}^{s}u_{i}||^{2}_{L^{2q}(\mathbb{R}^{n})}\Big)^{1/2}\\
&\leq cR^{\frac{n}{2}(1-\frac{1}{q})-\frac{1}{2}}||J^{s}f_{i}||_{L^{2q}(\mathbb{R}^{n})}\\
&\leq cR^{\frac{n}{2}(1-\frac{1}{q})-\frac{1}{2}}||f_{i}||_{L^{p}(\mathbb{R}^{n})}\leq c||f_{i}||_{L^{p}(\mathbb{R}^{n})},
\end{aligned}
\end{equation}
where $0\leq\frac{1}{p}-\frac{1}{2q}\leq\frac{s}{n}$, $\frac{1}{q}\geq1-\frac{1}{n}$ and $i=1,2$. The estimate for $D_{x}^{i}u_{3}$ can follow the same line and argument as those given in \cite[Theorem 3.1]{Ruiz1993}, so that we conclude that
\begin{equation}
\Big(\frac{1}{R}\int_{B_{R}}|D_{x}^{s}u_{3}(x)|^{2}dx\Big)^{1/2}\leq c||f||_{p},
\end{equation}
for $p$ such that $\frac{1}{p}-\frac{1}{p'}\geq 1-(n-1)(\frac{1}{p}-\frac{1}{2})$. As a consequence, we obtain the estimate of (\ref{local nablau}). The proof of (\ref{local u}) follows the same lines as the proof of (\ref{local nablau}), it is suffices to use the fact that $u_{i}$ is pointwise bounded by the Bessel potential $J^{2s}f_{i}$ when $i=1,2$.
\end{proof}

\subsection{$L^{\infty}(\mathbb{R}^{n})$ estimate for resolvent of fractional Helmholtz operator}

\begin{proof} [\bf Proof of Theorem \ref{thm5}]
Denote by $K(x)$ the Schwartz kernel of the resolvent $[(-\Delta)^{s}-(1+i\varepsilon)]^{-1}$, it follows from the proof of Theorem 1.4. in \cite{Yao2016} that
\begin{equation}\label{K(x)}
|K(x)|\leq \left \{
\begin{aligned}
&C|x|^{2s-n}, ~&&0<|x|\leq 1,\\
&C|x|^{\frac{1-n}{2}}, ~&&|x|>1,
\end{aligned}
\right.
\end{equation}
and
\begin{equation*}
|\nabla K(x)|\leq \left \{
\begin{aligned}
&c|x|^{2s-1-n}, ~&&0<|x|\leq 1,\\
&c|x|^{\frac{1-n}{2}}, ~&&|x|>1.
\end{aligned}
\right.
\end{equation*}
Then, we claim that
\begin{equation}
|||K|\ast f||_{L^{\infty}_{\tau(\alpha)}(\mathbb{R}^{n})}\leq C||f||_{L^{\infty}_{\alpha}(\mathbb{R}^{n})},~~~|||\nabla K|\ast f||_{L^{\infty}_{\tau(\alpha)}(\mathbb{R}^{n})}\leq C||f||_{L^{\infty}_{\alpha}(\mathbb{R}^{n})}.
\end{equation}
Indeed, we can calculate that
\begin{equation}
\begin{aligned}
|(K\ast f)(x)|&\leq \int_{\mathbb{R}^{n}}||K(y)|||f(x-y)|dy\\
&\leq C||f||_{L^{\infty}_{\alpha}(\mathbb{R}^{n})}\Big(\int_{B_{1}(0)}|y|^{2s-n}\langle x-y\rangle^{-\alpha}dy+\int_{\mathbb{R}^{n}\setminus B_{1}(0)}|y|^{\frac{1-n}{2}}\langle x-y\rangle^{-\alpha}dy\Big).
\end{aligned}
\end{equation}
For $|x|\leq 4$, it is easy to see that
\begin{equation}\label{I0}
\begin{aligned}
|(|K|\ast f)(x)|&\leq C||f||_{L^{\infty}_{\alpha}(\mathbb{R}^{n})}\Big(\int_{B_{1}(0)}|y|^{2s-n}dy+\int_{\mathbb{R}^{n}\setminus B_{1}(0)}|y|^{\frac{1-n}{2}-\alpha}dy\Big)\\
&\leq C||f||_{L^{\infty}_{\alpha}(\mathbb{R}^{n})}\Big(\int_{0}^{1}|r|^{2s-1}dr+\int_{1}^{\infty}|r|^{\frac{1-n}{2}-\alpha+n-1}dr\Big)\\
&\leq C||f||_{L^{\infty}_{\alpha}(\mathbb{R}^{n})},
\end{aligned}
\end{equation}
where we use that fact that $1\leq s\leq\frac{n}{2}$ and $\alpha>\frac{n+1}{2}>\frac{n-1}{2}$. Similarly, we have
\begin{equation}\label{I 0}
\begin{aligned}
|(|\nabla K|\ast f)(x)|&\leq C||f||_{L^{\infty}_{\alpha}(\mathbb{R}^{n})}\Big(\int_{0}^{1}|r|^{2s-2}dr+\int_{1}^{\infty}|r|^{\frac{1-n}{2}-\alpha+n-1}dr\Big)\leq C||f||_{L^{\infty}_{\alpha}(\mathbb{R}^{n})}.
\end{aligned}
\end{equation}
For $|x|>4$, we have
\begin{equation}
I_{1}=\int_{B_{1}(0)}|y|^{2s-n}\langle x-y\rangle^{-\alpha}dy\leq C|x|^{-\alpha}.
\end{equation}
While for the estimate outside $B_{1}(0)$, it has been computation in the proof of \cite[Lemma 2.1]{Chen2021}, that is
\begin{equation}
\begin{aligned}
I_{2}:=&\int_{B_{\frac{|x|}{2}}(0)\setminus B_{1}(0)}|y|^{\frac{1-n}{2}}\langle x-y\rangle^{-\alpha}dy\leq C|x|^{-\alpha+\frac{n+1}{2}};\\
I_{3}:=&\int_{B_{\frac{|x|}{2}}(x)}|y|^{\frac{1-n}{2}}\langle x-y\rangle^{-\alpha}dy\leq C|x|^{-\tau(\alpha)};\\
I_{4}:=&\int_{\mathbb{R}^{n}\setminus\big(B_{\frac{|x|}{2}}(0)\cup B_{\frac{|x|}{2}}(x)\big)}|z|^{\frac{1-n}{2}}\langle x-y\rangle^{-\alpha}dy\leq C|x|^{-\alpha+\frac{n+1}{2}}.
\end{aligned}
\end{equation}
 Since $-\tau(\alpha)\geq\mathrm{max}\{-\frac{n-1}{2},-\alpha,-\alpha+\frac{n+1}{2}\}$, one may combine these estimates with (\ref{I0}) to obtain the first statement of our claim. Moreover, for $|x|>4$, we also have
\begin{equation}
\begin{aligned}
\widetilde{I}_{1}=\int_{B_{1}(0)}|y|^{2s-1-n}\langle x-y\rangle^{-\alpha}dy\leq C|x|^{-\alpha}\leq C\langle x\rangle^{-\alpha}.
\end{aligned}
\end{equation}
Therefore, we find by (\ref{I 0}) that
\begin{equation}
 |(|\nabla K|\ast f)(x)|\leq C||f||_{L^{\infty}_{\alpha}}(\widetilde{I}_{1}+I_{2}+I_{3}+I_{4})\leq C\langle x\rangle^{-\tau(\alpha)}||f||_{L^{\infty}_{\alpha}(\mathbb{R}^{n})},~~\mathrm{for~all}~x\in\mathbb{R}^{n}.
\end{equation}
This implies that the second statement of our claim also holds. What's more, follow the same line in the proof of Proposition 1.1 in \cite{Chen2021}, the resolvent operator $\mathcal{R}^{s}_{\lambda}$ is a compact linear map.
\end{proof}

\section{Complex valued solutions for fractional Helmholtz equation}\label{complex}
\begin{proof}[\bf Proof of Theorem \ref{thm3}]
Define the operator $\mathcal{R}^{s}_{\lambda,\varphi}$ by
\begin{equation}
\mathcal{R}^{s}_{\lambda,\varphi}(u)(x)=\varphi(x)+\mathcal{R}^{s}_{\lambda}(|u|^{t-1}u)(x).
\end{equation}
We now show that for an appropriate $a>0$ the map $\mathcal{R}^{s}_{\lambda,\varphi}:B_{a}(L^{q}(\mathbb{R}^{n}))\longrightarrow B_{a}(L^{q}(\mathbb{R}^{n}))$ is a contraction, where $q$ satisfies the assumptions in $(i)-(v)$.

From Theorem \ref{thm1}, for any exponent $q$ satisfies the assumptions in $(i)-(v)$, there exist $p=\frac{q}{t}$ such that
\begin{equation}
\begin{aligned}
||\mathcal{R}^{s}_{\lambda,\varphi}(u)||_{L^{q}(\mathbb{R}^{n})}&\leq ||\varphi||_{L^{q}(\mathbb{R}^{n})}+||\mathcal{R}^{s}_{\lambda}(|u|^{t-1}u)||_{L^{q}(\mathbb{R}^{n})}\\
&\leq C||\varphi||_{L^{q}(\mathbb{R}^{n})}+|||u|^{t-1}u||_{L^{p}(\mathbb{R}^{n})}\\
&\leq C||\varphi||_{L^{q}(\mathbb{R}^{n})}+||u||_{L^{q}(\mathbb{R}^{n})}^{t}\\
&\leq C(\varepsilon+a^{t})<a,~~~\forall u\in B_{a}(L^{q}(\mathbb{R}^{n})).
\end{aligned}
\end{equation}
Moreover, due to the linearity of operator $\mathcal{R}^{s}_{\lambda}$, and the estimates in
Theorem \ref{thm1}, with the same indices $p$ and $q$, we obtain the following chain of
inequalities
\begin{equation}
\begin{aligned}
||\mathcal{R}^{s}_{\lambda,\varphi}(u)-\mathcal{R}^{s}_{\lambda,\varphi}(w)||_{L^{q}(\mathbb{R}^{n})}
&=||\mathcal{R}^{s}(|u|^{t-1}u-|w|^{t-1}w)||_{L^{q}(\mathbb{R}^{n})}\\
&\leq C|||u|^{t-1}u-|w|^{t-1}w||_{L^{p}(\mathbb{R}^{n})}\\
&\leq C||(|u|^{t-1}+|w|^{t-1})|u-w|||_{L^{p}(\mathbb{R}^{n})}\\
&\leq C(||u||_{L^{q}(\mathbb{R}^{n})}^{t-1}+||w||_{L^{q}(\mathbb{R}^{n})}^{t-1})||u-w||_{L^{q}(\mathbb{R}^{n})}\\
&\leq Ca^{t-1}||u-w||_{L^{q}(\mathbb{R}^{n})}<||u-w||_{L^{q}(\mathbb{R}^{n})},
\end{aligned}
\end{equation}
where we use the H\"{o}lder's inequality in obtaining the third inequality.

Therefore, the map $\mathcal{R}^{s}_{\lambda,\varphi}$ is a contraction in $B_{a}(L^{q}(\mathbb{R}^{n}))$, and consequently, there exists a unique $u\in B_{a}(L^{4}(\mathbb{R}^{n}))$ which satisfies $u=\mathcal{R}^{s}_{\lambda,\varphi}(|u|^{t-1}u)+\varphi$, that is a solution of (\ref{main system 1}).
\end{proof}

\begin{proof}[\bf Proof of Theorem \ref{thm6}]
Let $\varphi\in X:=L^{\infty}(\mathbb{R}^{n})$. We write (\ref{main system 1}) as a fixed point equation \begin{equation}
u=\mathcal{A}(u)~~\mathrm{in}~X
\end{equation}
with the nonlinear operator
\begin{equation}
\mathcal{A}:X\longrightarrow X,~~\mathcal{A}[w]=\mathcal{R}^{s}_{\lambda}(N_{f}(w))+\varphi,
\end{equation}
where $N_{f}(u)(x)=f(x,u(x))$ is a superposition operator. Since $\alpha>\frac{n+1}{2}$, we may fix $\alpha'\in(\frac{n+1}{2},\alpha)$. By Lemma 3.1 in \cite{Chen2021}, the nonlinear operator $N_{f}:X\longrightarrow L^{\infty}_{\alpha'}(\mathbb{R}^{n})$ is well defined and continuous. Moreover, $\mathcal{R}^{s}_{\lambda}:L^{\infty}_{\alpha'}(\mathbb{R}^{n})\longrightarrow X$ is compact by Theorem \ref{thm5}. Consequently, $\mathcal{A}$ is a compact and continuous operator. Moreover, let $\mathcal{F}\subset L^{\infty}(\mathbb{R}^{n})$ be  a set of function $u$ which solve the equation
\begin{equation}
u=\mu(\mathcal{R}^{s}_{\lambda}N_{f}(u)+\varphi)~~\mathrm{for~some~}\mu\in[0,1].
\end{equation}
Then, we easily deduce that $\mathcal{F}$ is bounded in $L^{\infty}(\mathbb{R}^{n})$, see the same analysis in the proof of Proposition 5.1 in \cite{Chen2021}. Hence a Schaefer's fixed point theorem  \cite[Chapter 9.2.2]{Evans1991} implies that $\mathcal{A}$ has a fixed point.

Furthermore, the Lipschitz condition (\ref{lips}) implies that
\begin{equation}
||\mathcal{A}(u)-\mathcal{A}(v)||_{X}=||\mathcal{R}^{s}_{\lambda}(N_{f}(u)-N_{f}(v))||\leq \kappa_{\alpha}||N_{f}(u)-N_{f}(v)||_{L^{\infty}_{\alpha}}\leq \kappa_{\alpha} l_{\alpha}||u-v||_{X},
\end{equation}
with $\kappa_{\alpha}l_{\alpha}<1$. Hence $\mathcal{A}$ is a contraction, and thus it has a unique fixed point in $X$.
\end{proof}

\begin{proof}[\bf Proof of Theorem \ref{thm8}]
Let $X=L^{\infty}(\mathbb{R}^{n})$, consider the nonlinear operator $\mathcal{B}:X\longrightarrow X$, $\mathcal{B}(u):=u-\mathcal{R}^{s}_{\lambda}N_{f}(u)$. Since $N_{f}(0)=0$, then $\mathcal{B}(0)=0$. Moreover, since $f(x,u)$ is continuous and differentiable, then $N_{f}:X\longrightarrow L^{\infty}_{\alpha'}$ is also differentiable, see the specific proof in Lemma 3.2 in \cite{Chen2021}. Therefore, $\mathcal{B}$
is a differentiable map. Moreover, by the assumption on $f(x,u)$, we have
\begin{equation}
\mathcal{B}'(0)=id-\mathcal{R}^{s}_{\lambda}N'_{f}(0)=id\in\mathcal{L}_{\mathbb{R}}(X,X).
\end{equation}
Consequently, $\mathcal{B}$ is a diffeomorphism between open neighborhoods $U,V\subset X$ of zero, and this implies that our claim.
\end{proof}

\section{The nonvanishing property}\label{non}
As we introduced before, we use the dual variational methods rather than the contraction mapping argument to detect the real valued solutions for (\ref{main system 1}), some compactness problem will arise naturally, hence, it is necessary to establish a nonvanishing property for resolvent.
\begin{lem}\label{vanishing}
Let $n\geq3$, $\frac{n}{n+1}<s<\frac{n}{2}$ and $\frac{2(n+1)}{n-1}<p<\frac{2n}{n-2s}$. Moreover, let $(v_{n})_{n}\subset L^{p'}(\mathbb{R}^{n})$ be a bounded sequence satisfying $\mathop{\mathrm{lim~sup}}\limits_{n\longrightarrow\infty}|\int_{\mathbb{R}^{n}}v_{n}\mathcal{R}_{\lambda}^{s}v_{n}|dx>0$.
Then there exists $R>0$, $\zeta>0$ and a sequence $(x_{n})_{n}\subset \mathbb{R}^{n}$ such that, up to a subsequence,
\begin{equation}
\int_{B_{R}(x_{n})}|v_{n}|^{p'}dx\geq\zeta~~\mathrm{for~all}~n.
\end{equation}
\end{lem}
The proof is similar to the classical case in \cite[Section 3]{Evequoz2015}, and we use a cut-off skill that revive the method of Guti\'{e}rrez in \cite{Gutierrez2004}. Fix $\psi\in\mathcal{S}$ such that $\widehat{\psi}\in\mathcal{C}^{\infty}(\mathbb{R}^{n})$ is radial, $0\leq\widehat{\psi}\leq1$, $\widehat{\psi}(\xi)=1$ for $||\xi|-1|\leq\frac{1}{6}$ and $\widehat{\psi}(\xi)=0$ for $||\xi|-1|\geq\frac{1}{4}$. We then write $K=K_{1}+K_{2}$ with
\begin{equation}
K_{1}:=\psi\ast K,~~~K_{2}=K-K_{1}.
\end{equation}
From (\ref{K(x)}) it then follows, by making $C_{0}>0$ large necessary, that $K_{1}\in\mathcal{C}^{\infty}(\mathbb{R}^{n})$ and
\begin{equation}
|K_{1}(x)|\leq C_{0}(1+|x|)^{\frac{1-n}{2}}~~~~\mathrm{for}~x\in\mathbb{R}^{n}.
\end{equation}
This is particular implies that $|K_{2}(x)|=|[K-K_{1}](x)|\leq2C_{0}|x|^{2s-n}$ for $|x|\leq1$. Moreover, since $\widehat{K_{2}}=(1-\widehat{\psi})\widehat{K}$ and $\widehat{K}=(|\xi|^{2s}-1-i0)^{-1}$, we have $\widehat{K_{2}}\in\mathcal{C}^{\infty}(\mathbb{R}^{n})$ with $\widehat{K_{2}}(\xi)=(|\xi|^{2s}-1)^{-1}$ for $|\xi|\geq\frac{5}{4}$. This implies that $\partial^{\gamma}\widehat{K_{2}}\in L^{1}(\mathbb{R}^{n})$ for all $\gamma$ such that $|\gamma|>n-2s$, which gives $|K_{2}(x)|\leq \kappa_{\beta}|x|^{-\beta}$, $x\in\mathbb{R}^{n}$ for all $\beta>n-2s$ with some constant $\kappa_{\beta}>0$. In particular, by making $C_{0}>0$ large if necessary, we have
\begin{equation}\label{osci}
|K_{2}(x)|\leq C_{0}\mathrm{min}\{|x|^{2s-n},|x|^{-n}\}~~~\mathrm{for~}x\in\mathbb{R}^{n}\setminus\{0\}.
\end{equation}

\begin{proof}
For $2<p<\frac{2n}{n-2s}$, we claim that for any bounded sequence $(v_{n})_{n}\in L^{p'}(\mathbb{R}^{n})$ such that
\begin{equation}
\mathop{\mathrm{lim}}\limits_{n\longrightarrow\infty}\Big(\mathop{\mathrm{sup}}\limits_{y\in\mathbb{R}^{n}}
\int_{B_{\rho}(y)}|v_{n}|^{p'}dx\Big)=0~~~\mathrm{for~all~}\rho>0,
\end{equation}
we have
\begin{equation}
\int_{\mathbb{R}^{n}}v_{n}[K_{2}\ast v_{n}]dx\longrightarrow 0~~~\mathrm{as}~n\longrightarrow \infty.
\end{equation}
Indeed, setting $A_{R}:=\{x\in\mathbb{R}^{n}:\frac{1}{R}\leq|x|\leq R\}$ and $D_{R}:=\mathbb{R}^{n}\setminus A_{R}$ for $R>1$, we derive from (\ref{osci}) that
\begin{equation}
||K_{2}||_{L^{\frac{p}{2}}(D_{R})}\longrightarrow 0~~~\mathrm{as}~R\longrightarrow\infty,
\end{equation}
since $1<\frac{p}{2}<\frac{n}{n-2s}$. Hence, by Young's inequality,
\begin{equation}\label{K2}
\mathop{\mathrm{sup}}\limits_{n\in\mathbb{N}}\big|\int_{\mathbb{R}^{n}}v_{n}[(1_{D_{R}}K_{2})\ast v_{n}]dx\big|\leq ||K_{1}||_{L^{\frac{p}{2}}(D_{R})}\mathop{\mathrm{sup}}\limits_{x\in\mathbb{N}}||v_{n}||^{2}_{L^{p'}(\mathbb{R}^{n})}\longrightarrow 0~~\mathrm{as}~R\longrightarrow\infty.
\end{equation}
On the other hand, decomposing $\mathbb{R}^{n}$ into disjoint $N-$cubes $\{Q_{l}\}_{l\in\mathbb{N}}$ of side length $R$, and considering for each $l$ the $N-$cube $Q'_{l}$ with the same center as $Q_{l}$ but with side length $3R$, we find,
\begin{equation}
\begin{aligned}
\Big|\int_{\mathbb{R}^{n}}v_{n}[(1_{A_{R}}K_{2})\ast v_{n}]dx\Big|&\leq\mathop{\sum}\limits^{\infty}_{l=1}\int_{Q_{l}}\Big(\int_{\frac{1}{R}<|x-y|<R}|K_{2}(x-y)||v_{n}(x)||v_{n}(y)|dx \Big)dx\\
&\leq CR^{n-2s}\mathop{\sum}\limits^{\infty}_{l=1}\int_{Q_{l}}\int_{Q_{l}}\Big(\int_{Q_{l}}|v_{n}(x)||v_{n}(y)|dy\Big)dx\\
&\leq CR^{n-2s+\frac{2n}{p}}\mathop{\sum}\limits^{\infty}_{l=1}\int_{Q_{l}}\Big(\int_{Q_{l}}|v_{n}(x)|^{p'}dx\Big)^{\frac{2}{p'}}\\
&\leq CR^{n-2s+\frac{2n}{p}}\big[\mathop{\mathrm{sup}}\limits_{l\in\mathbb{N}}\int_{Q_{l}}|v_{n}|^{p'}dx\big]^{\frac{2}{p'}-1}
\mathop{\sum}\limits^{\infty}_{l=1}\int_{Q'_{l}}|v_{n}(x)|^{p'}dx\\
&\leq CR^{n-2s+\frac{2n}{p}}\big[\mathop{\mathrm{sup}}\limits_{y\in\mathbb{R}^{n}}\int_{y\in\mathbb{R}^{n}}|v_{n}|^{p'}dx\big]^{\frac{2}{p'}-1}
3^{N}||v_{n}||^{p'}_{p'}.\\
\end{aligned}
\end{equation}
Hence, we also have
\begin{equation}\label{K22}
\mathop{\mathrm{lim}}\limits_{n\longrightarrow\infty}\int_{\mathbb{R}^{n}}v_{n}[(1_{A_{R}}K_{2})\ast v_{n}]dx=0~~~\mathrm{for~every~}~R>0.
\end{equation}
Combine (\ref{K2}) and (\ref{K22}), our claim holds.

For $p>\frac{2(n+1)}{n-1}$, we also claim that for any bounded sequence $(v_{n})_{n}\in L^{p'}(\mathbb{R}^{n})$ such that
\begin{equation}
\mathop{\mathrm{lim}}\limits_{n\longrightarrow\infty}\Big(\mathop{\mathrm{sup}}\limits_{y\in\mathbb{R}^{n}}
\int_{B_{\rho}(y)}|v_{n}|^{p'}dx\Big)=0~~~\mathrm{for~all~}\rho>0,
\end{equation}
we have
\begin{equation}
\int_{\mathbb{R}^{n}}v_{n}[K_{1}\ast v_{n}]dx\longrightarrow 0~~~\mathrm{as}~n\longrightarrow \infty.
\end{equation}
Indeed, since $|K_{1}(x)|\leq C_{0}(1+|x|)^{\frac{1-n}{2}}$, this proof follows the same line in Proposition 3.3 and  Lemma 3.4 in \cite{Evequoz2015}.
\end{proof}

\section{Variational Setting}\label{var}
In order to set up the variational framework, we begin with some preliminary work.
\subsection{Regularity Lemma}
As we can see, the $L^{q}$ complex valued solutions can be obtained by a bounded estimate for resolvent. Actually, these solutions also have a local strong regularity.
\begin{lem}\label{regularity000}
Let $n\geq3$, $\frac{n}{n+1}<s<\frac{n}{2}$, $\frac{2(n+1)}{n-1}<p<\frac{2n}{n-2s}$ and let $f\in L^{p'}(\mathbb{R}^{n})$. Then $u=\mathcal{R}_{\lambda}^{s}f\in W^{2s,p'}_{\mathrm{loc}}(\mathbb{R}^{n})\cap L^{p}(\mathbb{R}^{n})$ is a strong solution of $(-\Delta)^{s}u-u=f$ in $\mathbb{R}^{n}$. Moreover, for every $r>0$, there exists a constant $C>0$ depending on $r,s,n,p$, such that for all $x_{0}\in\mathbb{R}^{n}$,
\begin{equation}\label{regularity0}
||u||_{W^{2s,p'}(B_{r}(x_{0}))}\leq C\big(||u||_{L^{p}(\mathbb{R}^{n})}+||f||_{L^{p'}(\mathbb{R}^{n})}  \big).
\end{equation}
Furthermore,$~~$\\
(i)~If $f\in L^{p'}(\mathbb{R}^{n})\cap L^{q}_{\mathrm{loc}}(\mathbb{R}^{n})$ and $u\in L^{q}_{\mathrm{loc}}(\mathbb{R}^{n})$ for some $q\in(1,\infty)$, then $u\in W^{2s,q}_{\mathrm{loc}}(\mathbb{R}^{n})$, and for every $r>0$, there exists a constant $D>0$ depending only on $r,n,p,q$ such that
\begin{equation}\label{regularity00}
||u||_{W^{2s,q}(B_{r}(x_{0}))}\leq D\big(||u||_{L^{q}(B_{2r}(x_{0}))}+||f||_{L^{q}(B_{2r}(x_{0}))} \big)
\end{equation}
for every $x_{0}\in\mathbb{R}^{n}$.$~~$\\
(ii)~If $f\in L^{p'}(\mathbb{R}^{n})\cap L^{q}(\mathbb{R}^{n})$ and $u\in L^{q}(\mathbb{R}^{n})$ for some $q\in(1,\infty)$, then $u\in W^{2s,q}(\mathbb{R}^{n})$.
\end{lem}
\begin{proof}
Firstly, we claim that, for any $f\in L^{p'}(\mathbb{R}^{n})$, the equation $(-\Delta)^{s}u-u=f$ holds in the sense of distributional. Indeed, assume that $f\in\mathcal{S}$, then $\mathcal{R}_{\lambda}^{s}f\in\mathcal{S}'$ is given by
\begin{equation}
\langle\mathcal{R}_{\lambda}^{s}f,\varphi\rangle=\mathop{\mathrm{lim}}\limits_{\varepsilon\longrightarrow0^{+} }\int_{\mathbb{R}^{n}}\varphi(x)\mathbb{F}^{-1}\big(\frac{\widehat{f}(\cdot)}{|\cdot|^{2s}-1-i\varepsilon}\big)dx=
\mathop{\mathrm{lim}}\limits_{\varepsilon\longrightarrow0^{+}}\int_{\mathbb{R}^{n}}\frac{\check{\varphi}\widehat{f}(\xi)}{|\xi|^{2s}-1-i\varepsilon}
d\xi
\end{equation}
for all $\varphi\in\mathcal{S}$, where $\check{\varphi}$ is an abbreviation for $\mathbb{F}^{-1}(\varphi)$. On the other hand, since $\big|\frac{i\varepsilon}{|\xi|^{2s}-1-i\varepsilon}\big|\leq1$ for every $\xi\in\mathbb{R}^{n}$ and $\varepsilon>0$, and $\mathop{\mathrm{lim}}\limits_{\varepsilon\longrightarrow0^{+}}\frac{i\varepsilon}{|\xi|^{2s}-1-i\varepsilon}=0$ for $\xi\in\mathbb{R}^{n}$ with $|\xi|\neq1$, it then follows from Lebesgue's Theorem that
\begin{equation}
\mathop{\mathrm{lim}}\limits_{\varepsilon\longrightarrow0^{+}}\int_{\mathbb{R}^{n}}\frac{i\varepsilon}{|\xi|^{2s}-1-i\varepsilon}
g(\xi)d\xi=0~~~~\mathrm{for~every}~g\in\mathcal{S}.
\end{equation}
Hence, setting $u=\mathcal{R}_{\lambda}^{s}f$, we obtain for every $\varphi\in\mathcal{S}$
\begin{equation}
\begin{aligned}
\langle(-\Delta)^{s}u-u,\varphi\rangle=\langle\mathcal{R}_{\lambda}^{s}f,(-\Delta)^{s}\varphi-\varphi\rangle
&=\mathop{\mathrm{lim}}\limits_{\varepsilon\longrightarrow0^{+}}\int_{\mathbb{R}^{n}}\frac{\widehat{f}(\xi)\check{\varphi}(\xi)(|\xi|^{2s}-1)}{|\xi|^{2s}-1-\i\varepsilon}d\xi\\
&=\mathop{\mathrm{lim}}\limits_{\varepsilon\longrightarrow0^{+}}\int_{\mathbb{R}^{n}}\frac{\widehat{f}(\xi)\check{\varphi}(\xi)(|\xi|^{2s}-1-\i\varepsilon)}{|\xi|^{2s}-1-i\varepsilon}d\xi
=\mathcal{R}^{s}_{\lambda}f(x)\varphi(x)dx=\langle f,\varphi\rangle.
\end{aligned}
\end{equation}
This implies that the equation $(-\Delta)^{s}u-u=f$ holds in the distributional sense for any $f\in\mathcal{S}$. Now let $f\in L^{p'}(\mathbb{R}^{n})$ and consider a sequence $(f_{n})_{n}\subset\mathcal{S}$ with $||f_{n}-f||_{L^{p'}(\mathbb{R}^{n})}\longrightarrow 0$ as $n\longrightarrow\infty$. Then $u_{n}=\mathcal{R}_{\lambda}^{s}f_{n}$ solves $(-\Delta)^{s}u_{n}-u_{n}=f_{n}$ in the distributional sense, and $u_{n}\longrightarrow u$ in $L^{p}(\mathbb{R}^{n})$ by Theorem \ref{thm1}. Consequently, $(-\Delta)^{s}u_{n}-u_{n}=\longrightarrow f$ and $u_{n}\longrightarrow u$ in $\mathcal{S}'$ as $n\longrightarrow \infty$, this proves our claim.

Secondly, taking $x_{0}\in\mathbb{R}^{n}$, $r>0$ and considering the mollification $(u_{\varepsilon})_{\varepsilon>0}$ of $u=\mathcal{R}_{\lambda}^{s}f$, i.e., $u_{\varepsilon}:=\rho_{\varepsilon}\ast u$ where $\rho_{\varepsilon}(x)=\varepsilon^{-N}\rho(\frac{x}{\varepsilon})$, $x\in\mathbb{R}^{n}$ for some function $\rho\in\mathcal{C}_{c}^{\infty}(\mathbb{R}^{n})$ satisfying $\rho(x)\geq0$ for all $x\in\mathbb{R}^{n}$, $\mathrm{supp}(\rho)\subset B_{1}$ and $\int_{\mathbb{R}^{n}}\rho dx=1$. Since $u\in L^{p}(\mathbb{R}^{n})$, we deduce that $u\in L^{p'}(B_{r}(x_{0}))$ and consequently, $u_{\varepsilon}\longrightarrow u$ in $L^{p'}(B_{r}(x_{0}))$ as $\varepsilon\longrightarrow 0^{+}$. Similarly, considering the mollification $(f_{\varepsilon})_{\varepsilon>0}$ of $f$, we see that $f_{\varepsilon}\longrightarrow f$ in $L^{p'}(\mathbb{R}^{n})$ and therefore also in $L^{p'}(B_{r}(x_{0}))$, as $\varepsilon\longrightarrow 0^{+}$. Based on the properties of the mollification of $L^{p}-$ functions and of tempered distributions, with respect to differential operators with constant coefficients (see [\cite{Rudin1991}]), we yield that
\begin{equation}
(-\Delta)^{s}u_{\varepsilon}-u_{\varepsilon}=(-\Delta)^{s}(u\ast\rho_{\varepsilon})-(u\ast\rho_{\varepsilon})=((-\Delta)^{s}u-u)\ast\rho_{\varepsilon}=f\ast\rho_{\varepsilon}
=f_{\varepsilon}~~\mathrm{in}~\mathbb{R}^{n}.
\end{equation}
Therefore, the local regularity theory (see \cite[Theorem 1.4]{Biccari2017}) for fractional elliptic equation shows that the existence, for all $r>0$, of some constant $C>0$, depending only on $r,p$ and $N$,  such that
\begin{equation}\label{regularity}
||u_{\varepsilon}||_{W^{2s,p'}(B_{r}(x_{0}))}\leq C\big(||u_{\varepsilon}||_{L^{p'}(B_{2r}(x_{0}))}+||f_{\varepsilon}||_{L^{p'}(B_{2r}(x_{0}))}\big)~~~\mathrm{for~all}~\varepsilon>0.
\end{equation}
Then choosing some sequence $(\varepsilon_{n})_{n}\subset(0,\infty)$ such that $\varepsilon_{n}\longrightarrow 0$ as $n\longrightarrow\infty$ and replacing $u_{\varepsilon}$ by $u_{\varepsilon_{n}}-u_{\varepsilon_{m}}$ in (\ref{regularity}) gives that $(u_{\varepsilon_{n}})_{n}$ is a Cauchy sequence in $W^{2s,p'}(B_{r}(x_{0}))$ and therefore, there exists $w\in W^{2,p'}(B_{r}(x_{0}))$ such that $u_{\varepsilon_{n}}\longrightarrow w$ in $W^{2,p'}(B_{r}(x_{0}))$ as $n\longrightarrow\infty$. This also implies that $u_{\varepsilon_{n}}\longrightarrow w$ in $L^{p'}(B_{r}(x_{0}))$, Moreover, since $w=u$ a.e. in $B_{r}(x_{0})$, it follows that $u\in W^{2s,p'}(B_{r}(x_{0}))$ and $u$ solves the equation $(-\Delta)^{s}u-u=f$ almost everywhere in $B_{r}(x_{0})$. Furthermore, (\ref{regularity}) gives
\begin{equation}
\begin{aligned}
||u||_{W^{2s,p'}(B_{r}(x_{0}))}
&\leq C\big(||u||_{L^{p'}(B_{2r}(x_{0}))}+||f||_{L^{p'}(B_{2r}(x_{0}))}\big)\\
&\leq \tilde{C}\big(||u||_{L^{p}(\mathbb{R}^{n})}+||f||_{L^{p'}(\mathbb{R}^{n})}\big),
\end{aligned}
\end{equation}
where $\tilde{C}=C\mathrm{max}\{1,[\omega_{n}(2r)^{n}]^{\frac{p-2}{p}}\}$ and $\omega_{n}$ denotes the volume of the unit ball in $\mathbb{R}^{n}$. Since $r>0$ and $x_{0}\in\mathbb{R}^{n}$ were arbitrarily chosen, it follows that $u\in W^{2s,p'}_{\mathrm{loc}}(\mathbb{R}^{n})$ is a strong solution of $(-\Delta)^{s}u-u=f$ and, for every $r>0$, there exists a constant $\widetilde{C}>0$ depending only on $r,n,p$ such that (\ref{regularity0}) holds for all $x_{0}\in\mathbb{R}^{n}$.

The proof of (ii) follows the same line in \cite[Proposition A.1]{Evequoz2015} and we omit it here. While the claim in (iii) follows a global fractional Calder\'{o}n-Zygmund estimate, which has been presented in Abdellaoui, Fernandez, Leonori and Younes \cite{Abdellaoui2021}.
\end{proof}

\subsection{Local Compactness for Resolvent}
Being interested in the real-valued solutions of (\ref{main system 1}), all function spaces are assumed to consist of real-valued functions and we write $\Psi_{\lambda}^{s}:=\mathrm{Re}~K(x)$ for the real part of the Schwartz kernel $K(x)$ of the resolvent. By Theorem \ref{thm1}, we know that the linear operator $\mathbf{R}:L^{p'}(\mathbb{R}^{n})\longrightarrow L^{p}(\mathbb{R}^{n})$, $\mathbf{R}(v):=\Psi_{\lambda}^{s}\ast v$ is bounded.  Setting $v=Q^{\frac{1}{p'}}|u|^{p-2}u$, we are led to consider the equation
\begin{equation}\label{refor1}
|v|^{p'-2}v=Q^{\frac{1}{p}}[\Psi^{s}_{\lambda}\ast(Q^{\frac{1}{p}}v)]~~~\mathrm{in}~\mathbb{R}^{n}.
\end{equation}
To set up a dual variational framework for (\ref{main system 1}), it is necessary for us to study the operator $Q^{\frac{1}{p}}[\Psi^{s}_{\lambda}\ast(Q^{\frac{1}{p}}v)]$ on the right-hand side of (\ref{refor1}). Note $\mathbf{K}_{p}=Q^{\frac{1}{p}}[\Psi^{s}_{\lambda}\ast(Q^{\frac{1}{p}}v)]$, we have the following local compactness lemma for $\mathbf{K}_{p}$.
\begin{lem}\label{compact}
Let $n\geq3$, $\frac{n}{n+1}<s<\frac{n}{2}$, $\frac{2(n+1)}{n-1}\leq p\leq\frac{2n}{n-2s}$ and consider $Q\in L^{\infty}(\mathbb{R}^{n})$, satisfying $Q(x)\geq 0$ for a.e. $x\in\mathbb{R}^{n}$. Then $\mathbf{K}_{p}$ is symmetric in the sense that $\int_{\mathbb{R}^{n}}w\mathbf{K}_{p}(v)=\int_{\mathbb{R}^{n}}v\mathbf{K}_{p}(w)$ for all $v,w\in L^{p'}(\mathbb{R}^{n})$. Moreover, if $p<\frac{2n}{n-2s}$, we have $~~$\\
(i)~For any bounded and measurable set $B\subset\mathbb{R}^{n}$, the operator $1_{B}\mathbf{K}_{p}$ is compact. Here $1_{B}$ denote the characteristic function of the set $B$.$~~$\\
(ii)~If in addition, $\mathop{\mathrm{ess~sup}}\limits_{|x|\geq R}Q(x)\longrightarrow 0$ as $R\longrightarrow\infty$, then $\mathbf{K}_{p}$ itself is compact.
\end{lem}
\begin{proof}
By the argument in Lemma \ref{regularity000}, we know that $u\in\mathrm{W}^{2s,p'}_{\mathrm{loc}}(\mathbb{R}^{n})\cap L^{p}(\mathbb{R}^{n})$ is a strong solution of $(-\Delta)^{s} u-\lambda u=f$ in $\mathbb{R}^{N}$ for any $f\in L^{p'}(\mathbb{R}^{n})$. Furthermore, since $p<\frac{2n}{n-2s}$, we then have the compact embedding $W^{2s,p'}(B)\hookrightarrow L^{\frac{np'}{n-2sp'}}(B)\hookrightarrow L^{p}(B)$. Plugging this fact into the proof of
\cite[Lemma 4.1]{Evequoz2015}, we deduce that the operator $1_{B}\mathbf{K}_{p}$ is compact. And the remains follows the lines of the proof of Lemma 4.1 in \cite{Evequoz2015}.
\end{proof}

\subsection{Mountain Pass Structure}
Now, let us consider the energy functional
\begin{equation}
\begin{aligned}
J(v)&=\frac{1}{p'}\int_{\mathbb{R}^{n}}|v|^{p'}dx-\frac{1}{2}\int_{\mathbb{R}^{n}}Q(x)^{\frac{1}{p}}v(x)\mathbf{R}_{\lambda}^{s}(Q^{\frac{1}{p}}v)(x)dx\\
&=\frac{1}{p'}||v||_{p'}^{p'}-\frac{1}{2}\int_{\mathbb{R}^{n}}v\mathbf{K}_{p}(v)dx
\end{aligned}
\end{equation}
for $v\in L^{p'}(\mathbb{R}^{n})$. By the symmetric properties of the operator $\mathbf{K}_{p}$, we easily deduce that $J\in\mathcal{C}^{1}(L^{p'}(\mathbb{R}^{n}),\mathbb{R})$ with
\begin{equation}
J'(v)w=\int_{\mathbb{R}^{n}}\big(|v|^{p'-2}v-\mathbf{K}_{p}(v)\big)w~dx~~~\mathrm{for~all~}v,w\in L^{p'}(\mathbb{R}^{n}).
\end{equation}
Moreover, the functional $J$ has the so-called mountain pass geometry.
\begin{lem}\label{mountain pass}
Let $n\geq3$, $\frac{n}{n+1}<s<\frac{n}{2}$ and $\frac{2(N+1)}{N-1}\leq p\leq\frac{2n}{n-2s}$, consider $Q\in L^{\infty}(\mathbb{R}^{n})$ and $Q\not\equiv0$, such that $Q(x)\geq0$ for a.e. $x\in\mathbb{R}^{n}$.$~~$\\
(i)~There exists $\delta>0$ and $0<\rho<1$ such that $J(v)\geq \delta>0$ for all $v\in L^{p'}(\mathbb{R}^{n})$ with $||v||_{p'}=\rho$.$~~$\\
(ii)~There is $v_{0}\in L^{p'}(\mathbb{R}^{n})$ such that $||v_{0}||_{L^{p'}(\mathbb{R}^{n})}>1$ and $J(v_{0})<0$.$~~$\\
(iii)~Every Palais-Smale sequence for $J$ is bounded in $L^{p'}(\mathbb{R}^{n})$.
\end{lem}
\begin{proof}
(i)~By the boundedness of operator $\mathbf{R}$, there exists some constant $C>0$ such that $||\mathbf{K}_{p}(v)||_{L^{p}(\mathbb{R}^{n})}\leq C||v||_{L^{p'}(\mathbb{R}^{n})}$ for all $v\in L^{p'}(\mathbb{R}^{n})$. Hence, if $||v||_{L^{p'}(\mathbb{R}^{n})}=\rho$, we obtain
\begin{equation}
\begin{aligned}
J(v)&=\frac{1}{p'}\rho^{p'}-\frac{1}{2}\int_{\mathbb{R}^{n}}v\mathbf{K}_{p}(v)dx\geq \frac{1}{p'}\rho^{p'}-\frac{1}{2}\big(||\mathbf{K}_{p}(v)||_{L^{p}(\mathbb{R}^{n})}||v||_{L^{p'}(\mathbb{R}^{n})}\\
&=\frac{1}{p'}\rho^{p'}-\frac{\rho}{2}||\mathbf{K}_{p}(v)||_{L^{p}(\mathbb{R}^{n})}\geq\frac{1}{p'}\rho^{p'}-\frac{C}{2}\rho^{2}>0,\\
\end{aligned}
\end{equation}
where we use the fact of $p'<2$ and we chose $\rho>0$ small enough. $~~$\\
(ii)~From the representation formulations in \cite{Yao2016}, it follows that there exists $r>0$ such that $\Psi_{\lambda}^{s}>0$ for all $x\in B_{2R}(0)$. Moreover, since $Q\geq 0$ a.e. on $\mathbb{R}^{n}$ and $Q\not\equiv0$, the metric density of the set $\omega_{Q}:=\{x\in\mathbb{R}^{n}:Q(x)\geq0\}$ is 1 for almost every point from this set. Therefore, there exists $x_{0}\in\mathbb{R}^{n}$ and $0<\rho<r$ such that $\omega_{Q}\cap B_{\frac{\rho}{2}(x_{0})}$ has positive measure. Choosing $z\in\mathcal{C}_{c}^{\infty}(\mathbb{R}^{n})$ with $\mathrm{supp}~z\subset B_{\rho}(x_{0})$, $0\leq z\leq 1$ in $\mathbb{R}^{n}$ and $z=1$ in $B_{\frac{\rho}{2}}(x_{0})$, the definition of $\mathbf{K}_{p}$ implies that
\begin{equation}
\begin{aligned}
\int_{\mathbb{R}^{n}}z\mathbf{K}_{p}zdx&=\int_{\mathbb{R}^{n}}\int_{\mathbb{R}^{n}}Q(x)^{\frac{1}{p}}z(x)\Psi_{\lambda}^{s}(x-y)Q(y)^{\frac{1}{p}}z(y)dydx\\
&\geq \int_{B_{\frac{\rho}{2}(x_{0})}}\int_{B_{\frac{\rho}{2}(x_{0})}}\Psi_{\lambda}^{s}(x-y)Q(x)^{\frac{1}{p}}Q(y)^{\frac{1}{p}}dxdy>0.
\end{aligned}
\end{equation}
Hence taking $t>0$ large enough, we obtain
\begin{equation}
\begin{aligned}
J(tz)&=\frac{t^{p'}}{p'}\int_{\mathbb{R}^{n}}|z|^{p'}dx-\frac{t^{2}}{2}\int_{\mathbb{R}^{n}}z\mathbf{K}_{p}zdx\\
&=t^{2}\Big(\frac{1}{p't^{2-p'}}\int_{\mathbb{R}^{n}}|z|^{p'}dx-\frac{1}{2}\int_{\mathbb{R}^{n}}z\mathbf{K}_{p}zdx\Big)<0.
\end{aligned}
\end{equation}
(iii)~For every Palais-Smale sequence $(v_{n})_{n}\subset L^{p'}(\mathbb{R}^{n})$, we have
\begin{equation}
\begin{aligned}
+\infty&>\mathop{\mathrm{sup}}\limits_{n}|J(v_{n})|\geq J(v_{n})=(\frac{1}{p'}-\frac{1}{2})||v_{n}||_{p'}^{p'}+\frac{1}{2}J'(v_{n})v_{n}\\
&\geq(\frac{1}{p'}-\frac{1}{2})||v_{n}||_{p'}^{p'}-\frac{1}{2}||J'(v_{n})||_{L^{p'}(\mathbb{R}^{n})^{\ast}}||v_{n}||_{L^{p'}(\mathbb{R}^{n})}\\
&\geq(\frac{1}{p'}-\frac{1}{2})||v_{n}||_{p'}^{p'}~~\mathrm{as}~n\longrightarrow\infty.
\end{aligned}
\end{equation}
This implies $(v_{n})_{n}$ is bounded since $1<p'<2$.
\end{proof}

\section{Real solutions for fractional Helmholtz equation}\label{real}
\subsection{Existence of solutions in the compact case}
For the case $Q(x)\longrightarrow 0$ as $|x|\longrightarrow \infty$, we shall prove the existence of infinitely many pairs $\{\pm u\}$ of critical points for $J$ using a variant of the symmetric Mountain Pass Theorem, see \cite{Ambrosetti1973}. Therefore, we need more properties of $\mathbf{K}_{p}$ and the functional $J$.
\begin{lem}
For every $m\in\mathbb{N}$, there exist an $m-$dimendsional subspace $\mathcal{W}\subset\mathcal{C}_{c}^{\infty}(\mathbb{R}^{n})$ with the following properties:$~~$\\
(i)~$\int_{\mathbb{R}^{n}}v\mathbf{K}_{p}vdx>0$ for all $v\in\mathcal{W}\setminus\{0\}$.$~~$\\
(ii)~There exists $R=R(\mathcal{W})>0$ such that $J(v)\leq 0$ for every $v\in\mathcal{W}$ with $||v||_{L^{p'}(\mathbb{R}^{n})}\geq R$.
\end{lem}

\begin{proof}
Since $Q\not\equiv0$, there exists a point $x_{0}$ of density one for the set $\{Q>0\}$. Without loss of generality, we may assume that $x_{0}=0$. Then for $\delta>0$ sufficiently small we have
\begin{equation}\label{0set}
|Q^{-1}(0)\cap B_{\delta}(0)|\leq(\frac{1}{4m^{2}})^{n}|B_{\delta}(0)|.
\end{equation}
Let
\begin{equation}
\Psi_{\lambda}^{s,\ast}(\tau):=\mathop{\mathrm{inf}}\limits_{B_{\tau}(0)\setminus\{0\}}\Psi_{\lambda}^{s}~~~\mathrm{and}
\Psi_{\lambda,\ast}^{s}(\tau):=||\Psi_{\lambda}^{s}||_{L^{\infty}(\mathbb{R}^{n}\setminus B_{\tau}(0))}~~\mathrm{for}~~\tau>0.
\end{equation}
Since $\Psi_{\lambda}^{s}$ is bounded outside of every neighborhood of zero and $\Psi_{\lambda}^{s}(x)|x|^{n-2s}$ tends to a positive constant as $|x|\longrightarrow 0$ by the representation of fractional resolvent in \cite{Yao2016}, we may fix $\delta>0$ such that (\ref{0set}) holds and that
\begin{equation}\label{ast}
\Psi_{\lambda}^{s,\ast}(\tau)>(m-1)\Psi_{\lambda,\ast}^{s}(m\tau)~~~~\mathrm{for}~\tau\in(0,\delta].
\end{equation}
Moreover, it is easy to see that there exists $m$ disjoint open balls $B^{1},...,B^{m}\subset B_{\delta}(0)$ of diameter $\tau:=\frac{\delta}{m^{2}}$ such that
\begin{equation}
\mathrm{dist}(B^{i},B^{j}):=\mathrm{inf}\{|x-y|:x\in B^{i},y\in B^{i}\}\geq\frac{\delta}{m}.
\end{equation}
Since $|B^{i}|=(\frac{1}{2m^{2}})^{n}|B_{\delta}(0)|$ for $i=1,...,m$, then by (\ref{0set}), we also have
\begin{equation}\label{0set1}
|B^{i}\cap\{Q>0\}|>0~~~\mathrm{for}~i=1,...,m.
\end{equation}
Now, fix functions $z_{i}\in \mathcal{C}_{c}^{\infty}(\mathbb{R}^{n})$, $i=1,...,m$ such that $z_{i}>0$ in $B^{i}$ and $z_{i}\equiv 0$ in $\mathbb{R}^{n}\setminus B^{i}$. Moreover, we let $\mathcal{W}$ denote the span of $z_{1},...,z_{m}$. Then any $v\in\mathcal{W}\setminus\{0\}$ can be written as $v=\mathop{\sum}\limits_{i=1}^{m}a_{i}z_{i}$ with $a=(a_{1},...,a_{m})\in\mathbb{R}^{m}\setminus\{0\}$, and thus we have
\begin{equation}
\begin{aligned}
\int_{\mathbb{R}^{n}}v\mathbf{K}_{p}vdx&=\mathop{\sum}\limits_{i,j=1}^{m} a_{i}a_{j}\int_{B^{i}}\int_{B^{j}}\Psi_{\lambda}^{s}Q(x)^{\frac{1}{p}}Q(y)^{\frac{1}{p}}z_{i}(x)z_{j}(y)dxdy\\
&\geq \Psi_{\lambda}^{s,\ast}(\tau)\mathop{\sum}\limits_{i=1}^{m}a_{i}^{2}\Big(\int_{B^{i}}Q(x)^{\frac{1}{p}}z_{i}(x)dx\Big)^{2}\\
&~~~-\Psi_{\lambda,\ast}^{s}(m\tau)\mathop{\sum}\limits_{i,j=1,\i\neq j}^{m}|a_{i}||a_{j}|\Big(\int_{B^{i}}Q(x)^{\frac{1}{p}}z_{i}(x)dx\Big)\Big(\int_{B^{j}}Q(x)^{\frac{1}{p}}z_{j}(x)dx\Big)\\
&\geq \Psi_{\lambda}^{s,\ast}(\tau)\mathop{\sum}\limits_{i=1}^{m}a_{i}^{2}\Big(\int_{B^{i}}Q(x)^{\frac{1}{p}}z_{i}(x)dx\Big)^{2}\\
&~~~-\frac{\Psi_{\lambda,\ast}^{s}(m\tau)}{2}\mathop{\sum}\limits_{i,j=1,i\neq j}\Big[a_{i}^{2}\Big(\int_{B^{i}}Q(x)^{\frac{1}{p}}z_{i}(x)\Big)^{2}+a_{j}^{2}\Big(\int_{B^{j}}Q(x)^{\frac{1}{p}}z_{j}(x)\Big)^{2}\Big]\\
=&\mathop{\sum}\limits_{i=1}^{m}\Big(\Psi^{s,\ast}_{\lambda}(\tau)-(m-1)\Psi^{s}_{\lambda,\ast}(m\tau)\Big)a_{i}^{2}\Big(\int_{B^{i}}Q(x)^{\frac{1}{p}}z_{i}(x)dx\Big)^{2}>0,
\end{aligned}
\end{equation}
as a consequence of (\ref{ast}) and (\ref{0set1}). This proves (i).

By the continuity of $\mathbf{K}_{p}$, we have
\begin{equation}
m_{\mathcal{W}}:=\mathop{\mathrm{inf}}\limits_{v\in\mathcal{W},||v||_{L^{p}(\mathbb{R}^{n})}=1}\int_{\mathbb{R}^{n}}v\mathbf{K}_{p}vdx>0.
\end{equation}
Therefore, we obtain that
\begin{equation}
J(v)=\frac{||v||_{L^{p'}(\mathbb{R}^{n})}^{p'}}{p'}-\frac{1}{2}\int_{\mathbb{R}^{n}}v\mathbf{K}_{p}vdx\leq ||v||_{L^{p'}(\mathbb{R}^{n})}^{p'}(\frac{1}{p'}-\frac{1}{2}||v||_{L^{p'}(\mathbb{R}^{n})}^{2-p'}m_{\mathcal{W}})~~~\mathrm{for}~v\in\mathcal{W}.
\end{equation}
Taking $R:=\Big(\frac{2}{m_{\mathcal{W}}p'}\Big)^{\frac{1}{2-p'}}$, we have $J(v)\leq 0$ for every $v\in\mathcal{W}$.
\end{proof}

\begin{lem}\label{lem10}
Let $n\geq 3$, $\frac{n}{n+1}<s<\frac{n}{2}$, $\frac{2(n+1)}{n-1}<p<\frac{2n}{n-2}$, and let $Q\in L^{\infty}(\mathbb{R}^{n})$, $Q\geq0$, $Q\not\equiv0$ satisfy $\mathop{\mathrm{lim}}\limits_{|x|\longrightarrow\infty}Q(x)=0$. Then problem $u=\mathrm{Re}(\mathcal{R}_{\lambda}^{s}(Q(x)|u|^{p-2}u))$ admits a sequence of pairs $\pm u_{n}$ of solutions such that $||u_{n}||_{L^{p}(\mathbb{R}^{n})}\longrightarrow \infty$ as $n\longrightarrow \infty$.
\end{lem}

\begin{proof}
Let $(v_{n})_{n}\subset L^{p'}(\mathbb{R}^{n})$ be a Palais-Smale sequence, then by Lemma \ref{mountain pass}~(iii), we know that $(v_{n})_{n}$ is bounded in $L^{p'}(\mathbb{R}^{n})$, Hence, up to a subsequence, we may assume $v_{n}\rightharpoonup v\in L^{p'}(\mathbb{R}^{n})$, this also implies that $||v||_{L^{p'}(\mathbb{R}^{n})}\leq \mathop{\mathrm{lim~inf}}\limits_{n\longrightarrow \infty}||v_{n}||_{L^{p'}(\mathbb{R}^{n})}$. On the other hand, we easily obtain that
\begin{equation}
\begin{aligned}
\frac{1}{p'}||v||_{L^{p'}(\mathbb{R}^{n})}^{p'}&-\frac{1}{p'}||v_{n}||_{L^{p'}(\mathbb{R}^{n})}^{p'}\geq \int_{\mathbb{R}^{n}}|v_{n}|^{p'-2}v_{n}(v-v_{n})dx\\
&=J'(v_{n})(v-v_{n})+\int_{\mathbb{R}^{n}}v_{n}\mathbf{K}_{p}(v-v_{n})dx\longrightarrow 0,~~n\longrightarrow \infty,
\end{aligned}
\end{equation}
where we use the convexity of the function $t\longrightarrow |t|^{p'}$. As a consequence, we have
$||v||_{L^{p'}(\mathbb{R}^{n})}\leq \mathop{\mathrm{lim}}\limits_{n\longrightarrow \infty}||v_{n}||_{L^{p'}(\mathbb{R}^{n})}$, this implies that $v_{n}\longrightarrow v$ strongly in $L^{p'}(\mathbb{R}^{n})$, which means that $J(v)$ satisfies the Palais-Smale condition. Combining the previous lemma with the symmetric Mountain Pass Theorem, we then obtain the existence of nontrivial pairs $\{\pm v_{n}\}$ of critical point of $J$ with $J(v_{n})\longrightarrow \infty$ and thus $||v_{n}||_{L^{p'}(\mathbb{R}^{n})}\longrightarrow\infty$ as $n\longrightarrow\infty$. Remembering the setting $u_{n}:=\mathbf{R}_{\lambda}^{s}(Q^{\frac{1}{p}}v_{n})$, we have $v_{n}=Q^{\frac{1}{p'}}|u_{n}|^{p-2}u_{n}$, thus $||u_{n}||_{L^{p}(\mathbb{R}^{n})}\longrightarrow\infty$ as $n\longrightarrow\infty$.
\end{proof}

\subsection{Existence of solutions in the periodic case}
As we have showed in Lemma \ref{mountain pass}, the functional $J(v)$ satisfies the mountain pass geometry. Hence, we may define a mountain-pass level for $J(v)$ by
\begin{equation}
c:=\mathop{\mathrm{inf~max}}\limits_{\gamma\in\Gamma, t\in[0,1]}J(\gamma(t))
\end{equation}
where $\Gamma=\{\gamma\in C([0,1],L^{p'}(\mathbb{R}^{n})):\gamma(0)=0~\mathrm{and}~J(\gamma(1))<0\}$. Apparently, $\Gamma\neq\emptyset$ and $c>0$. To show that $c$ is a critical level of $J(v)$, we consider some Palais-Smale sequences $(v_{n})_{n}\subset L^{p'}(\mathbb{R}^{n})$ for $J(v)$ at level $c$, which can be easily found via a deformation Lemma. Moreover, these Palais-Smale sequences have been proved to be bounded in Lemma \ref{mountain pass} (iii). A long as we prove these sequences having a strong convergence subsequences in $L^{p'}(\mathbb{R}^{n})$, we then obtain the main conclusion.
\begin{lem}\label{lem11}
Let $n\geq 3$, $\frac{n}{n+1}<s<\frac{n}{2}$, $\frac{2(n+1)}{n-1}<p<\frac{2n}{n-2}$, Consider a nonnegative function $Q\in L^{\infty}(\mathbb{R}^{n})$, $Q\not\equiv0$ which is $\mathbb{Z}^{n}-$periodic on $\mathbb{R}^{n}$. Then problem $u=\mathrm{Re}(\mathcal{R}_{\lambda}^{s}(Q(x)|u|^{p-2}u))$ has a nontrivial solution $u\in L^{p}(\mathbb{R}^{n})$.
\end{lem}
\begin{proof}
Let $(v_{n})_{n}$ be a bounded Palais-Smale sequence for $J(v)$ at level $c$, that is $J(v_{n})\longrightarrow c$ and $J'(v_{n})v_{n}\longrightarrow 0$ as $n\longrightarrow \infty$. Then we easily deduce that
\begin{equation}
\mathop{\mathrm{lim}}\limits_{n\longrightarrow\infty}\int_{\mathbb{R}^{n}}Q^{\frac{1}{p}}v_{n}\mathbf{R}_{\lambda}^{s}(Q^{\frac{1}{p}}v_{n})dx
=\frac{2p'}{2-p'}\mathop{\mathrm{lim}}\limits_{n\longrightarrow \infty}[J(v_{n})-\frac{1}{p'}J'(v_{n})v_{n}]=\frac{2p'}{2-p'}c>0.
\end{equation}
Moreover, since $Q\in L^{\infty}(\mathbb{R}^{n})$, the sequence $(Q^{\frac{1}{p}}v_{n})_{n}$ is also bounded. Therefore, by the vanishing lemma \ref{vanishing}, there exists $R,\zeta>0$ and a sequence $(x_{n})_{n}\subset \mathbb{R}^{n}$ such that
\begin{equation}
\int_{B_{R}(x_{0})}|v_{n}|^{p'}dx\geq\zeta~~~~\mathrm{for~all}~n.
\end{equation}
By the periodicity of $Q$, this problem is invariance under translation. Hence we may set $w_{n}(x)=v_{n}(x+x_{n})$ for $x\in\mathbb{R}^{n}$, where $(w_{n})_{n}\subset L^{p'}(\mathbb{R}^{n})$ is a bounded sequence such that $J(w_{n})=J(v_{n})$ and $||J'(w_{n})||=||J'(v_{n})||$.
Moreover, up to a subsequence, we have $w_{n}\rightharpoonup w$ in $L^{p'}(\mathbb{R}^{n})$. Next, we claim that
\begin{equation}
1_{B_{R'}}|w_{n}|^{p'-2}w_{n}\longrightarrow 1_{B_{R'}}|w|^{p'-2}w~~\mathrm{strongly~in~}L^{p}(\mathbb{R}')~for ~every~R'>0.
\end{equation}
Indeed, fix $\varphi\in\mathcal{C}_{c}^{\infty}(B_{R'})\subset\mathcal{C}_{c}^{\infty}(\mathbb{R}^{n})$, then for $n,m\in\mathbb{N}$ we have
\begin{equation}
\begin{aligned}
\Big|&\int_{\mathbb{R}^{n}}\Big(|w_{n}|^{p'-2}w_{b}-|w_{m}|^{p'-2}w_{m}\Big)\varphi dx\Big|\\
&=\Big|J'(w_{n})\varphi-J'(w_{m})\varphi+\int_{B_{R'}}\varphi\mathbf{K}_{p}(w_{n}-w_{m})dx\Big|\\
&\leq ||J'(w_{n})-J"(w_{m})||_{\mathcal{L}^{p}(L^{p}(\mathbb{R}^{n}),\mathbb{R})}||\varphi||_{L^{p'}(\mathbb{R}^{n})}+
||1_{B_{R'}}\mathbf{K}_{p}(w_{n}-w_{m})||_{\mathcal{L}^{p}(L^{p}(\mathbb{R}^{n}),\mathbb{R})}||\varphi||_{L^{p'}(\mathbb{R}^{n})}.
\end{aligned}
\end{equation}
Since $\mathcal{C}_{c}^{\infty}(B_{R'})\subset L^{p'}(B_{R'})$ is dense, $||J'(w_{n})||\longrightarrow 0$ as $n\longrightarrow \infty$ and since $1_{B_{R'}}\mathbf{K}_{p}$ is a compact operator, we deduce that $|w_{n}|^{p'-2}w_{n}$ is a Cauchy sequence in $L^{p}(B_{R'})$, so that $|w_{n}|^{p'-2}w_{n}\longrightarrow \widetilde{w}$ strongly in $L^{p}(B_{R'})$ for some $\widetilde{w}\in L^{p}(B_{R'})$. Up to a subsequence, $|w_{n}|^{p'-2}w_{n}\longrightarrow \widetilde{w}$ and, equivalently, $w_{n}\longrightarrow |\widetilde{w}|^{p-2}\widetilde{w}$ a.e. on $B_{R'}$. By the uniqueness of the weak limit, it follows that $w=|\widetilde{w}|^{p-2}\widetilde{w}$, i.e. $\widetilde{w}=|w|^{p'-2}w$ on $B_{R'}$. This proves our claim.

As a consequence,
\begin{equation}
0<\zeta\leq\int_{B_{R}(x_{n})}=\int_{B_{R}}|w_{n}|^{p'}dx\longrightarrow\int_{B_{R}}|w|^{p'}dx~~\mathrm{as}~n\longrightarrow\infty,
\end{equation}
which implies $w\neq0$. We are going to show that $w$ is a critical point of $J(w)$. For every $\varphi\in\mathcal{C}_{c}^{\infty}(\mathbb{R}^{n})$, we have
\begin{equation}
\int_{\mathbb{R}^{n}}|w_{n}|^{p'-2}w_{n}\varphi dx\longrightarrow\int_{\mathbb{R}^{n}}|w|^{p'-2}w\varphi dx~~~\mathrm{as}~n\longrightarrow\infty.
\end{equation}
On the other hand since $\mathbf{K}_{p}$ is a bounded operator, we have
\begin{equation}
\int_{\mathbb{R}^{n}}\varphi\mathbf{K}_{p}(w_{n})dx\longrightarrow\int_{\mathbb{R}^{n}}\varphi\mathbf{K}_{p}(w)dx~~\mathrm{as}~n\longrightarrow\infty.
\end{equation}
Consequently,
\begin{equation}
\begin{aligned}
J'(w)\varphi&=\int_{\mathbb{R}^{n}}|w|^{p'-2}w\varphi dx-\int_{\mathbb{R}^{n}}\varphi\mathbf{K}_{p}(w)dx\\
&=\mathop{\mathrm{lim}}\limits_{n\longrightarrow\infty}\Big(\int_{\mathbb{R}^{n}}|w_{n}|^{p'-2}w_{n}\varphi-\int_{\mathbb{R}^{n}}\varphi \mathbf{K}_{p}(w_{n})dx\Big)=\mathop{\mathrm{lim}}\limits_{n\longrightarrow\infty}J'(w_{n})\varphi=0.
\end{aligned}
\end{equation}
Therefore, $w\in L^{P'}(\mathbb{R}^{n})$ is a nontrivial critical point of $J(w)$, so is $v$. Remembering the setting $u:=\mathbf{R}_{\lambda}^{s}(Q^{\frac{1}{p}}v)$, we have $v=Q^{\frac{1}{p'}}|u|^{p-2}u$, this implies that $u\in L^{p}(\mathbb{R}^{n})$ is a nontrivial solution for $u=\mathrm{Re}(\mathcal{R}_{\lambda}^{s}(Q(x)|u|^{p-2}u))$.
\end{proof}

\subsection{Strong solutions for the Helmholtz equation}
As we have showed that that problem $u=\mathrm{Re}(\mathcal{R}_{\lambda}^{s}(Q(x)|u|^{p-2}u))$ has a nontrivial weak solution $u\in L^{p}(\mathbb{R}^{n})$. In the following, we study the regularity of $u$.

\begin{lem}\label{strong}
Let $n\geq 3$, $\frac{n}{n+1}<s<1$, $\frac{2(n+1)}{n-1}<p<\frac{2n}{n-2s}$. Let $Q\in L^{\infty}(\mathbb{R}^{n})$, and consider a solution $u\in L^{p}(\mathbb{R}^{n})$ of $u=\mathrm{Re}(\mathcal{R}_{\lambda}^{s}(Q(x)|u|^{p-2}u))$. Then $u\in W^{2s,q}(\mathbb{R}^{n})\cap\mathcal{C}^{1,\alpha}(\mathbb{R}^{n})$ for all $p\leq q<\infty$, $0<\alpha<1$, and it is a strong solution of (\ref{main system 1}).
\end{lem}
\begin{proof}
Firstly, we shall use the Moser iteration technique to show that the solutions are bounded in $L^{\infty}(\mathbb{R}^{n})$. Indeed, since $Q\in L^{\infty}(\mathbb{R}^{n})$ and $\frac{2(n+1)}{n-1}<p<\frac{2n}{n-2s}$, then by Lemma \ref{regularity000}, it follows that $u\in W^{2s,p'}_{\mathrm{loc}}(\mathbb{R}^{n})$ and for every $x_{0}\in\mathbb{R}^{n}$,
\begin{equation}
||u||_{W^{2s,p'}(B_{2r}(x_{0}))}\leq \widetilde{C}\Big(||u||_{L^{p}(\mathbb{R}^{n})}+||Q||_{L^{\infty}(\mathbb{R}^{n})}||u||^{p-1}_{L^{p}(\mathbb{R}^{n})} \Big)
\end{equation}
with some constant $\widetilde{C}>0$, independent of $x_{0}$. Moreover, $u$ is a strong solution of (\ref{main system 1}). Using Sobolev's embedding theorem with the property $p'\geq\frac{2n}{n+2s}$, we obtain that $u\in W^{s,2}_{\mathrm{loc}}(\mathbb{R}^{n})$ with
\begin{equation}
||u||_{W^{s,2}(B_{2r}(x_{0}))}\leq\kappa\widetilde{C}\Big(||u||_{L^{p}(\mathbb{R}^{n})}+||Q||_{L^{\infty}(\mathbb{R}^{n})}||u||_{L^{p-1}(\mathbb{R}^{n})}  \Big)~~\mathrm{for~all}~x_{0}\in\mathbb{R}^{n},
\end{equation}
where the constant $\kappa$ is independent of $x_{0}$.

Consider now $L>0$, $\beta>1$ and a cut off function $\eta\in\mathcal{C}_{c}^{\infty}(\mathbb{R}^{N})$ with $\mathrm{supp}~\eta\subset B_{r}(x_{0})$, define $u_{L}=\mathrm{min}\{u,L\}$ and $\gamma(u)=\eta uu_{L}^{2(\beta-1)}$. Apparently, $\gamma$ is a increasing function and we have
\begin{equation}
(a-b)(\gamma(a)-\gamma(b))\geq 0~~~~\mathrm{for~~any}~~a,b\in\mathbb{R}.
\end{equation}
Define the functions
\begin{equation}
\Lambda(t)=\frac{|t|^{2}}{2}~~~\mathrm{and}~~~\Gamma(t)=\int_{0}^{t}(\gamma'(\tau))^{\frac{1}{2}}d\tau.
\end{equation}
Fix $a,b\in\mathbb{R}$ such that $a>b$. Then, from the above definitions and applying Jensen inequality we
get
\begin{equation}
\begin{aligned}
\Lambda'(a-b)(\gamma(a)-\gamma(b))&=(a-b)(\gamma(a)-\gamma(b))=(a-b)\int_{a}^{b}\gamma'(t)dt\\
&=(a-b)\int_{a}^{b}(\Gamma'(t))^{2}dt\geq \Big(\int_{a}^{b}\Gamma'(t)dt\Big)^{2}=\Gamma(a)-\Gamma(b).
\end{aligned}
\end{equation}
In similar fashion, we can prove that the above inequality is true for any $a\leq b$. Thus we can infer
that
\begin{equation}
\Lambda'(a-b)(\gamma(a)-\gamma(b))\geq |\Gamma(a)-\Gamma(b)|^{2},
\end{equation}
In particular, it follows that
\begin{equation}\label{gamma}
\begin{aligned}
|\Gamma(\eta u(x))-\Gamma(\eta u(y))|^{2}\leq |\eta u(x)-\eta u(y)|((\eta uu_{L}^{2(\beta-1)})(x)-(\eta uu_{L}^{2(\beta-1)})(y)).
\end{aligned}
\end{equation}
Therefore, taking $\gamma(u)=\eta uu_{L}^{2(\beta-1)}$ as test-function in (\ref{main system 1}), in view of (\ref{gamma}) we have
\begin{equation}
\begin{aligned}
&||\Gamma(\eta u)||_{W^{s,2}(B_{r}(x_{0}))}^{2}-\int_{B_{r}(x_{0})}\eta u^{2}u_{L}^{2(\beta-1)}dx\\
&~~\leq \int\int_{B_{r}(x_{0})\times B_{r}(x_{0})}\frac{u(x)-u(y)}{|x-y|^{n+2s}}((\eta uu_{L}^{2(\beta-1)})(x)-(\eta uu_{L}^{2(\beta-1)})(y))dxdy
-\int_{B_{r}(x_{0})}\eta u^{2}u_{L}^{2(\beta-1)}dx\\
&~~=\int_{B_{r}(x_{0})}Q(x)\eta |u|^{p}|u_{L}|^{2(\beta-1)}dx.
\end{aligned}
\end{equation}
Since $\Gamma(\eta u)\geq \frac{1}{\beta}\eta uu_{L}^{\beta-1}$, from the Sobolev inequality we can deduce that
\begin{equation}
||\Gamma(\eta u)||_{W^{s,2}(B_{r}(x_{0}))}^{2}\geq S_{\ast}|\Gamma(\eta u)|_{L^{2^{\ast}_{s}}(B_{r}(x_{0}))}^{2}\geq (\frac{1}{\beta})^{2}S_{\ast}|\eta uu_{L}^{\beta-1}|_{L^{2^{\ast}_{s}}(B_{r}(x_{0}))}^{2},
\end{equation}
where $S_{\ast}$ is the best constant of the fractional Sobolev embedding inequality. By the fact $Q(x)\in L^{\infty}(\mathbb{R}^{n})$ and $p>2$, we have
\begin{equation}
\begin{aligned}
|\eta uu_{L}^{\beta-1}|_{L^{2^{\ast}_{s}}(B_{r}(x_{0}))}^{2}&\leq C(||Q||_{L^{\infty}(\mathbb{R}^{n})})S_{\ast}^{-1}\beta^{2}\int_{B_{r}(x_{0})}\eta|u|^{p}|u_{L}|^{2(\beta-1)}dx\\
&\leq C_{1}\beta^{2}\Big(\int_{B_{r}(x_{0})}|u|^{2^{\ast}_{s}}\Big)^{\frac{p-2}{2^{\ast}_{s}}}
\Big(\int_{B_{r}(x_{0})}|\eta uu_{L}|^{\frac{22^{\ast}_{s}}{2^{\ast}_{s}-(p-2)}}\Big)^{\frac{2^{\ast}_{s}-(p-2)}{2^{\ast}_{s}}} .
\end{aligned}
\end{equation}
Setting $w_{L}=\eta uu_{L}^{\beta-1}$ and $\alpha^{\ast}=\frac{22^{\ast}_{s}}{2^{\ast}-(p-2)}$, we then have
\begin{equation}
||w_{L}||_{L^{2^{\ast}_{s}}(B_{r}(x_{0}))}^{2}\leq C_{2}\beta^{2}||w_{L}||_{L^{\alpha^{\ast}}(B_{r}(x_{0}))}^{2}.
\end{equation}
Now, we observe that if $u^{\beta}\in L^{\alpha^{\ast}}(B_{r}(x_{0}))$, from the definition of $w_{L}$, $u_{L}\leq u$, we obtain
\begin{equation}
||w_{L}||_{L^{2^{\ast}_{s}}(B_{r}(x_{0}))}^{2}\leq C_{3}\beta^{2}\big(\int_{B_{r}(x_{0})}u^{\beta\alpha^{\ast}}dx\big)^{\frac{2}{\alpha^{\ast}}}<\infty.
\end{equation}
Since $\eta\in\mathcal{C}_{c}^{\infty}(\mathbb{R}^{n})$ was chosen arbitrarily with $\mathrm{supp}~\eta\subset B_{r}(x_{0})$, then by passing to the limit as $L\longrightarrow\infty$, the Fatou's Lemma yields
\begin{equation}\label{step0}
||u||_{L^{\beta2^{\ast}_{s}}(B_{r}(x_{0}))}\leq C_{3}^{\frac{1}{2\beta}}\beta^{\frac{1}{\beta}}||u||_{L^{\beta\alpha^{\ast}}(B_{r}(x_{0}))}.
\end{equation}
Now, we set $\beta=\frac{2^{\ast}_{s}}{\alpha^{\ast}}>1$ and we observe that, being $u\in L^{2^{\ast}_{s}}(B_{r}(x_{0}))$, the above inequality holds for this choice of $\beta$. Then, by using the fact that $\beta^{2}\alpha^{\ast}=\beta2^{\ast}_{s}$, it follows that (\ref{step0}) holds with $\beta$ replaced by $\beta^{2}$. Therefore, we can see that
\begin{equation}
\begin{aligned}
||u||_{L^{\beta^{2}2^{\ast}_{s}}(B_{r}(x_{0}))}&\leq C_{3}^{\frac{1}{2\beta^{2}}}\beta^{\frac{2}{\beta^{2}}}||u||_{L^{\beta^{2}\alpha^{\ast}}(B_{r}(x_{0}))}\\
&\leq C_{3}^{\frac{1}{2}(\frac{1}{\beta}+\frac{1}{\beta^{2}})}\beta^{\frac{1}{\beta}+\frac{1}{\beta^{2}}}||u||_{L^{\beta\alpha^{\ast}}(B_{r}(x_{0}))}.
\end{aligned}
\end{equation}
Iterating this process, and recalling that $\beta\alpha^{\ast}:=2^{\ast}_{s}$, we can infer that for every $m\in\mathbb{N}$
\begin{equation}
\begin{aligned}
||u||_{L^{\beta^{m}2^{\ast}_{s}}(B_{r}(x_{0}))}&\leq C_{3}^{\sum^{m}_{j-1}\frac{1}{2\beta^{j}}}\beta^{\sum^{m}_{j-1}j\beta^{j-1}}||u||_{L^{2^{\ast}_{s}}(B_{r}(x_{0}))}.
\end{aligned}
\end{equation}
Taking the limit as $m\longrightarrow\infty$ we get
\begin{equation}
||u||_{L^{\infty}(B_{r}(x_{0}))}\leq C_{4}^{\gamma_{1}}\beta^{\gamma_{2}}<\infty,
\end{equation}
where $\gamma_{1}=\frac{1}{2}\sum^{\infty}_{j=1}\frac{1}{\beta^{j}}<\infty$ and $\gamma_{2}=\sum^{\infty}_{j=1}\frac{j}{\beta}<\infty$. This implies that $u\in L^{\infty}(B_{r}(x_{0}))$ for all $x_{0}\in\mathbb{R}^{n}$ with $\mathop{\mathrm{sup}}\limits_{x_{0}\in\mathbb{R}^{n}}||u||_{L^{\infty}(B_{\frac{1}{2}}(x_{0}))}<\infty$, i.e. $u\in L^{\infty}(\mathbb{R}^{n})$, as claimed. Applying Lemma \ref{regularity000} (ii) we then find that $u\in W^{2s,q}(\mathbb{R}^{n})$ for every $p\leq q<\infty$.
\end{proof}

\begin{proof}[\bf Proof of Theorem \ref{thm10}]
By Lemma \ref{lem10} and Lemma \ref{strong}, we obtain the strong solutions for (\ref{main system 1}).
\end{proof}

\begin{proof}[\bf Proof of Theorem \ref{thm11}]
The proof of Theorem \ref{thm11} easily follows from Lemma \ref{lem11} and Lemma \ref{strong}.
\end{proof}

\section*{Acknowledgements}
The research of Zifei Shen was partially supported by NSFC(12071438).

\section*{Declarations}
\textbf{Conflict of interest} The authors declare that they have no conflict of interest.

\end{document}